\documentclass[a4paper,12pt,leqno]{article}
\usepackage{soul}
\usepackage{graphicx, graphics}
\usepackage{syntonly}
\usepackage{amssymb, amsmath, amsfonts}
\usepackage{subfigure}
\usepackage{multirow}
\usepackage{color}
\usepackage{hyperref}

\usepackage{bm}
\newcommand{\bX}{\bm x}
\newcommand{\bu}{\bm u}
\newcommand{\bv}{\bm v}

\newcommand{\ud}{\, d}
\newcommand{\br}{\bm r}
\providecommand{\pfrac}[2]{\frac{\partial #1}{\partial #2}}
\providecommand{\tnabla}{\tilde{\nabla}}
\providecommand{\tpartial}{\tilde{\partial}}
\providecommand{\cF}{\mathcal F}
\providecommand{\fF}{\mathcal F}

\providecommand{\cS}{\mathcal S}

\providecommand{\bp}{\bm p}

\begin{document}
%\begin{frontmatter}
\title{Nonlocal strong forms of thin plate, gradient elasticity, magneto-electro-elasticity and phase field fracture by nonlocal operator method}
\author{Huilong Ren,Xiaoying Zhuang,\\
Erkan Oterkus,HeHua Zhu,Timon Rabczuk\\
Institute of Structural Mechanics,\\
 Bauhaus-University Weimar,
 99423 Weimar, Germany\\
 College of Civil Engineering,Tongji University,\\
  Shanghai 200092, China\\
Department of Naval Architecture, Ocean and Marine Engineering,\\
 PeriDynamics Research Centre,University of Strathclyde,\\
100 Montrose Street, Glasgow G4 0LZ, UK 
\thanks{huilong.ren@uni-weimar.de; zhuang@ikm.uni-hannover.de;timon.rabczuk@tdt.edu.vn}}
\maketitle
\begin{abstract}
The derivation of nonlocal strong forms for many physical problems remains cumbersome in traditional methods. In this paper, we apply the variational principle/weighted residual method based on nonlocal operator method for the derivation of nonlocal forms for elasticity, thin plate, gradient elasticity, electro-magneto-elasticity and phase field fracture method. The nonlocal governing equations are expressed as integral form on support and dual-support. The first example shows that the nonlocal elasticity has the same form as dual-horizon non-ordinary state-based peridynamics. The derivation is simple and general and it can convert efficiently many local physical models into their corresponding nonlocal forms. In addition, a criterion based on the instability of the nonlocal gradient is proposed for the fracture modelling in linear elasticity. Several numerical examples are presented to validate nonlocal elasticity and the nonlocal thin plate .
\end{abstract}
%\begin{keyword}
%energy form\sep weak form\sep variational principle\sep peridynamics\sep dual-support\sep fracture \sep explicit time integration
%\end{keyword}
%\end{frontmatter}

%Xueweng Wang \thanks{Email:xwwang@whut.edu.cn}
%\\State Key Laboratory of Advanced Technology\\
% for Materials Synthesis and Processing,\\ International School of Materials Science and Engineering,\\ Wuhan University of Technology, Wuhan, 430070, China\\

%\providecommand{\bX}{\mathbb X}
\providecommand{\bY}{\mathbb Y}
\providecommand{\bZ}{\mathbb Z}
\providecommand{\bI}{\mathbb I}
%Variational derivation of nonlocal thin plate and nonlocal gradient elasticity based on nonlocal operator method
\graphicspath{{ThinPlate/}{plate2d/}{Mode3frac/}}
\section{Introduction}
%Local theory
%

Classical continuum mechanics has achieved great success in describing the macro-scale properties of solid material based on the continuous medium hypothesis that the material is a continuous mass rather than as discrete particles. The assumption indicates that the substance of the object completely fills the space it occupies, without considering the inherent micro-structure of the material. Such a continuous medium hypothesis is not always valid in solid medium. Over the years, researchers found that many phenomena, such as size effect \cite{Toupin1962Jan}, length scale effect \cite{Mindlin1962micro}, skin/edge effect \cite{davydov2013molecular}, can not be well predicted by traditional continuum mechanics. These phenomena may be attributed to the nonlocal effect in the solid. In contrast with local theory whose mathematical language is partial differential derivatives defined at an infinitesimal point, nonlocal theory is formulated as integral form in a domain. %In micro-scale, the nonlocal effect grows more and more significant. On one side, the micro-structure of the medium can�t be described by conventional elasticity. %The micro-structure %This and so on, which fall into the category of nonlocal theory of mechanics.
%Another concern is the inhomogeneity of material in micro-scale, which violates the continuum assumption.

Classical continuum mechanics is regarded as a local theory. For solid mediums of multiple materials with material interface or discontinuity such as fracture, the partial differential operator is no longer well defined. Around the fracture front tip, the stress singularity happens for local theory. In order to model fracture and its evolution, various local theories have been proposed, for example, finite element method (FEM) \cite{areias2019finite}, extended finite element method \cite{sukumar2000extended}, phase-field fracture method \cite{nguyen2018modeling,ren2019explicit,zhou2020phase}, cracking particle method \cite{rabczuk2004cracking,zhang2018cracking},
extended finite element method \cite{majidi2018use}, numerical manifold method \cite{yang2019stability},
extended isogeometric analysis (XIGA) for three-dimensional crack \cite{singh2018bezier}, meshfree methods \cite{belytschko1994element,liu1995reproducing,huerta2018meshfree}. Another approach for fracture modeling is the nonlocal method. Compared with continuum mechanics without length scale, nonlocal theory takes into account the length scale explicitly and it is less sensitive to the inhomogeneity/discontinuity encountered in the materials due to its integral form.

Two general theories to account for the length scale of solid material, are the gradient elasticity \cite{mindlin1968first,Toupin1962Jan,Yang2002May,polizzotto2012gradient} and the nonlocal elasticity \cite{eringen1972nonlocal2,eringen1972nonlocal,Eringen1983Sep,eringen2012microcontinuum}. The gradient elasticity theory can be traced back to Cosserat theory in 1909 \cite{cosserat1909theorie}. It incorporates the length scale and higher order derivative of the displacement field. A variety of gradient elasticity theories have been proposed such as Mindlin solid theory \cite{mindlin1968first,Mindlin1962micro}, couple stress theory \cite{Toupin1962Jan,Toupin1964Jan}, modified couple stress \cite{Yang2002May,tsiatas2009new} and second-grade materials \cite{polizzotto2012gradient}.
In nonlocal elasticity, the stress tensor is based on the integral of the ``local'' stress field in a domain, in contrast with the local elasticity defining the stress based on the strain field at a point. Under certain circumstances, the nonlocal elasticity can be transformed into gradient elasticity \cite{Eringen1983Sep,dell2015origins}. 

%The idea of nonlocality has been extended in many directions. Yan et al proposed a Updated Lagrangian Particle Hydrodynamics based on nonlocal theory and applied it to model multiphase flows \cite{yan2020higher}

Among various nonlocal elasticity theories, Peridynamics (PD) \cite{Silling2000,Silling2007} have attracted the attention of the researchers in the fracture mechanic field. PD is based on the integral form well defined in domain with/without discontinuity. This salient feature enables PD a versatile method for fracture modeling \cite{foster2011energy,liu2018modeling,Zhou2016Oct,zhou2016numerical}. The origin of PD is the bond-based PD (BB-PD) with the Poisson ratio restriction. BB-PD can model 2D elasticity with Poisson ratio of 1/3 and 3D elasticity with Poisson ratio of 1/4. Many efforts have been dedicated to overcome this restriction, for example, PD with shear deformation \cite{ren2016new}, bond-rotation effect by \cite{zhu2017peridynamic}, PD with micropolar deformation \cite{diana2019bond}. The further development of PD is the state-based PD \cite{Silling2007,silling2010peridynamic}. Several treatments are developed to overcome the instability issue in non-ordinate state-based PD (NOSBPD), including, bond-associated higher-order stabilized model \cite{Gu2019Dec}, higher-order approximation \cite{yaghoobi2017higher}, stabilized non-ordinary state-based PD \cite{Silling2017Aug,li2018stabilized}, sub-horizon scheme \cite{chowdhury2019modified} and stress point method \cite{cui2020higher}.

In the spirit of nonlocality, PD has been extended in many directions, for example, dual-horizon PD \cite{Ren2015,ren2017dual}, peridynamic plate/shell theory \cite{taylor2015two,chowdhury2016peridynamic,dorduncu2019stress,zhang2020peridynamic},
mixed peridynamic Petrov-Galerkin method for compressible and incompressible hyperelastic material \cite{bode2020peridynamic,bode2020mixed},
phase field based peridynamic damage model for composite structures \cite{roy2017phase}, wave dispersion analysis of PD \cite{butt2017wave}, damage mechanism in PD \cite{yu2020energy}, coupling scheme for state-based PD and FEM \cite{bie2018coupling,d2019review}, higher-order peridynamic material models for elasticity \cite{chen2020higher}, to list a few.

%The development of dual-horizon peridynamics. Nonlocality is not a constant parameter for material.

Dual-horizon PD overcomes the restriction of constant horizon in PD, without introducing side effects for variable support size. Dual-horizon peridynamic formulation can be derived from the Euler-Lagrange equations \cite{wang2020derivation}. Based on the concept in nonlocal theory, we developed the Nonlocal Operator Method (NOM) as the generalization of dual-horizon PD. %Nonlocal operator methods, brief introduction, features, applications, capabilities.
NOM uses the nonlocal operators of integral form to replace the local partial differential operators of different orders. There are three versions of NOM, first-order particle-based NOM \cite{ren2020nonlocal,rabczuk2019nonlocal}, higher order particle-based NOM \cite{ren2020higher} and higher order NOM based on numerical integration \cite{ren2020nonlocalni}. {The particle-based version can be viewed as a special case of NOM with numerical integration when nodal integration is employed}. The nonlocal operators can be viewed as an alternative to the partial derivatives of shape functions in FEM. Combined with a variational principle or weighted residual method, NOM obtains the residual vector and tangent stiffness matrix in the same way as in FEM. NOM has been applied to the solutions of the Poisson equation in high dimensional space, von-Karman thin plate equations, fracture problems based on phase field \cite{ren2020higher}, waveguide problem in electromagnetic field \cite{rabczuk2019nonlocal}, gradient solid problem \cite{ren2020nonlocalni} and Cahn-Hilliard equation \cite{ren2020nonlocalCH}.

Although much progress in nonlocal methods has been achieved in the above mentioned literatures, the derivations for many physical problems remain cumbersome and complicated, see for example \cite{chowdhury2016peridynamic,Wang2018Mar,chen2020higher,javili2020computational}. In local theory, the local differential operator is a fundamental element for describing physical problems. In analogy, the nonlocal operators would be very beneficial for developing nonlocal theoretical models. The power of NOM in deriving nonlocal models remains largely unexplored. In addition, NOM based on implicit algorithms is relatively complicated in implementation and in this paper, we explore the explicit algorithm in solving the nonlocal models. Furthermore, we propose an instability criterion of the nonlocal gradient operator for the purpose of fracture modeling. The remaining of the paper is outlined as follows. In section 2, the second-order NOM in 2D/3D is formulated in detail. In section 3, we apply the NOM scheme combined with variational principle/weighted residual method to derive the nonlocal governing equations for elasticity, thin plate, gradient elasticity, electro-magneto-elasticity and phase field fracture model. The correspondence between local form and nonlocal form for higher order problems is discussed. In section 4, an instability criterion of nonlocal gradient is presented in the fracture modeling of linear elastic solid. The implementation of nonlocal solid and nonlocal thin plate is discussed in section 5. Several numerical examples for solid and thin plate are used to demonstrate the accuracy and efficiency of the current method in section 6. Last but not the least, some concluding remarks are presented.

\section{Second-order nonlocal operator method}
NOM uses the integral form to replace the partial differential derivatives of different orders. Although NOM can solve higher order linear/nonlinear problems in 2D/3D, we restrict our discussion in second-order NOM, which is sufficient for the nonlocal derivation of the physical problems to be studied in section 3.
\subsection{Support and dual-support}\label{sec:nom}
\begin{figure}[htp]
\centering
\subfigure[]{
\label{fig:Coord}
\includegraphics[width=.4\textwidth]{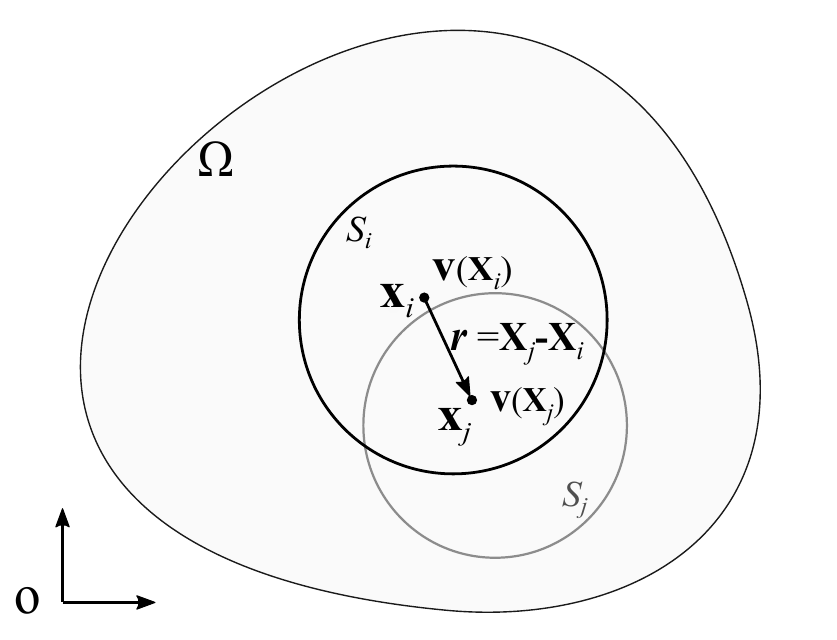}}
\vspace{.1in}
\subfigure[]{
\label{fig:4support}
\includegraphics[width=.35\textwidth]{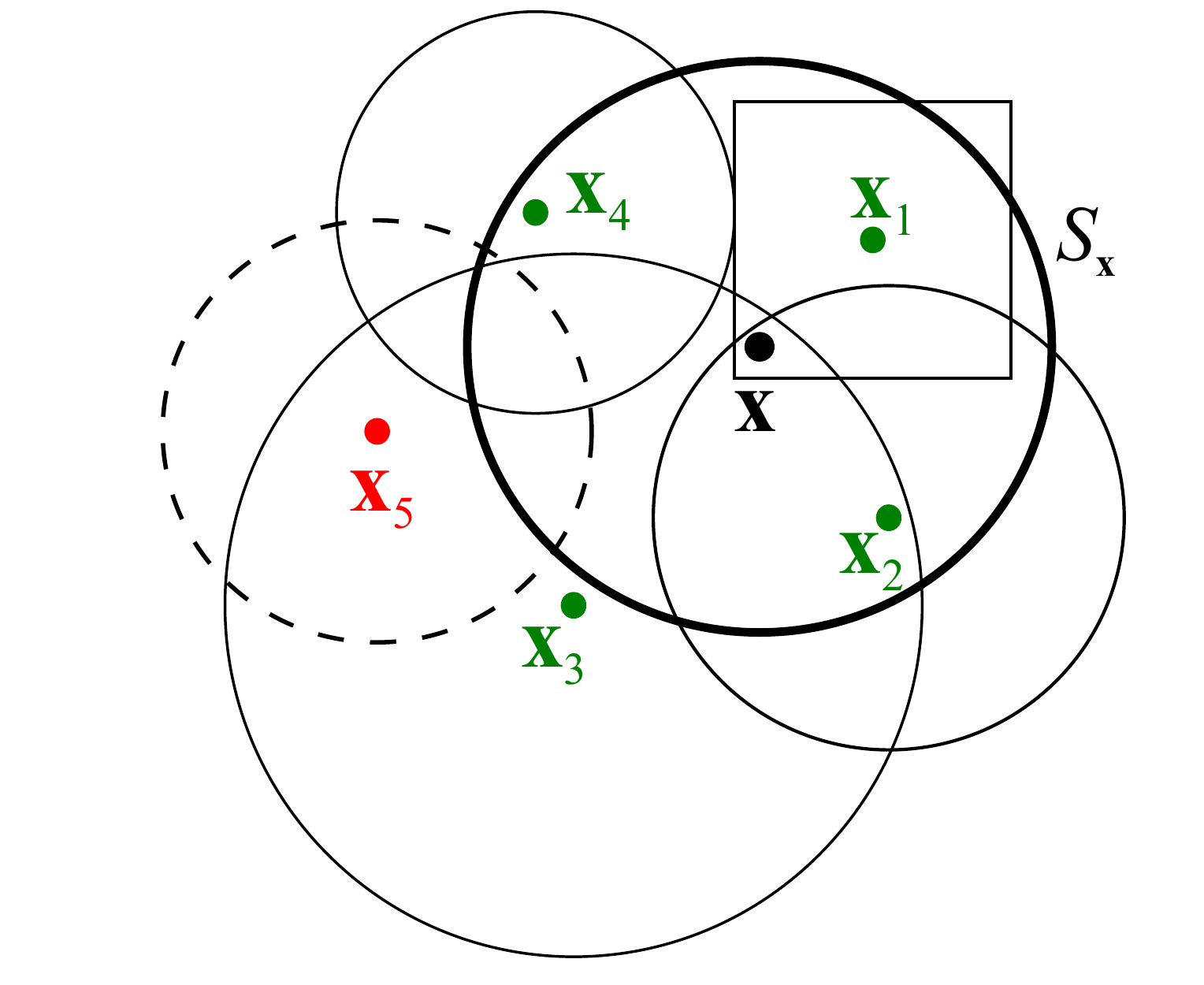}}\\
\vspace{.3in}
\caption{(a) Domain and notation. (b) Schematic diagram for support and dual-support, all shapes above are supports, $\cS_{\bX}=\{\bX_1,\bX_2,\bX_4\} $, $\cS_{\bX}'=\{\bX_1,\bX_2,\bX_3\}$.}
\end{figure}
Consider a domain as shown in Fig.\ref{fig:Coord}, let $\bX_i$ be spatial coordinates in the domain $\bm \Omega$; $\br_{ij}:=\bX_j-\bX_i$ is a spatial vector starting from $\bX_i$ to $\bX_j$; $\bv_i:=\bv(\bX_i,t)$ and $\bv_j:=\bv(\bX_j,t)$ are the field values for $\bX_i$ and $\bX_j$, respectively; $\bv_{ij}:=\bv_j-\bv_i$ is the relative field vector for spatial vector $\br$.

\textbf{Support} $\cS_{i}$ of point $\bX_i$ is the neighbourhood of point $\bX_i$. A point $\bX_j$ in support $\cS_i$ forms the spatial vector $\br(=\bX_j-\bX_i)$. The support in the NOM can be a spherical domain, a cube, semi-spherical domain and so on.

\textbf{Dual-support} is defined as a union of points whose supports include $\bX_i$, denoted by
\begin{align}
\cS_{i}'=\{\bX_j|\bX_i \in \cS_{j}\} \label{eq:dualsupport}.
\end{align}
Point $\bX_j$ forms the dual-vector $\br_{ji}(=\bX_i-\bX_j=-\br_{ij})$ in $\cS_{i}'$. On the other hand, $\br_{ji}$ is the spatial vector formed in $\cS_{j}$. It is worth mentioning that the size of the support of each point can be different. When the support sizes for all material points are the same, the dual-support is equal to the support. On the other hand, if the size of support varies for each point, the shape of dual-support can be quite irregular, even discontinuous for two adjacent points. One example to illustrate the support and dual-support is shown in Fig.\ref{fig:4support}.

\subsection{Dual property of dual-support}
For point $j\in \cS_i$, let $f_{ij}$ be a physical quantity, work conjugate to field difference $(u_j-u_i)$, the dual property of dual-support is
\begin{align}
\int_{\Omega} \int_{\cS_i} f_{ij} (u_j-u_i) \ud V_j\ud V_i=\int_{\Omega} \Big(\int_{\cS_i'} f_{ji} \ud V_j-\int_{\cS_i} f_{ij} \ud V_j\Big) u_i \ud V_i \label{eq:dualprop}
\end{align}

\noindent\textbf{Proof}:

Let the domain $\Omega$ be divided into $N$ non-overlapping particles, so that $\Omega=\sum_{i=1}^N \Delta V_i$, where $\Delta V_i$ is the volume assigned to particle $i$. Herein, $N$ can be arbitrarily large so that the $\Delta V_i$ is infinitesimal and the double summations of discrete form converge to the double integrals in continuous form.
\begin{align}
&\int_{\Omega} \int_{\cS_i} f_{ij} (u_j-u_i) \ud V_j\ud V_i\notag\\
\approx&\sum_{\Delta V_{i} \in \Omega}\sum_{\Delta V_j \in \cS_i}f_{ij} (u_j-u_i)\Delta V_j\Delta V_{i}\notag\\
=&\sum_{\Delta V_{i} \in \Omega}\sum_{ \Delta V_j \in \cS_i} f_{ij} u_j\Delta V_j\Delta V_{i}-\sum_{\Delta V_{i} \in \Omega}\sum_{ \Delta V_j \in \cS_i} f_{ij} u_i\Delta V_j\Delta V_{i}\notag\\
=&\sum_{\Delta V_{i} \in \Omega}\sum_{\Delta V_j \in \cS_i'} f_{ji} u_i\Delta V_j\Delta V_{i}-\sum_{\Delta V_{i} \in \Omega}\sum_{\Delta V_j \in \cS_i} f_{ij} u_i\Delta V_j\Delta V_{i}\notag\\
\approx&\int_{\Omega} \Big(\int_{\cS_i'} f_{ji} \ud V_j-\int_{\cS_i} f_{ij} \ud V_j\Big) u_i \ud V_i
\end{align}
In the third step, the dual-support is considered as follows. The term $f_{ij}$ with $u_j$ is the physical quantity from $i$'s support, but is added to particle $j$; since $j\in \cS_{i}$, $i$ belongs to the dual-support $\cS'_{j}$ of $j$; all terms $f_{ji}$ with $u_i$ are collected from any material point $j$ whose support contains $i$ and hence form the dual-support of $i$. Therefore, the dual property of the dual-support is proved.

When all points have the same size of support domains, i.e. $j \in \cS_i \leftrightarrow i\in \cS_j$, we have $\cS_i=\cS_i'$ for any point $i$ and then the dual property of dual-support by Eq.\ref{eq:dualprop} becomes
\begin{align}
\int_{\Omega} \int_{\cS_i} f_{ij} (u_j-u_i) \ud V_j\ud V_i=\int_{\Omega} \int_{\cS_i}( f_{ji}- f_{ij} ) u_i\ud V_j \ud V_i \label{eq:dualprop2}
\end{align}

Above equation is widely used in the derivation of nonlocal strong form from weak form. Such expression is valid in the continuum form as well as in discrete form. The dual property of dual-support is also proved in the dual-horizon peridynamics \cite{ren2017dual}. A simple example with $N=4$ to illustrate this property is given in \ref{sec:app1}.

\subsection{Nonlocal gradient and Hessian operator}
The local gradient operator and Hessian operator for a scalar-valued function $u$ have the forms in 2D
\begin{align}
\nabla u=\Big(u_{,x},u_{,y}\Big)^T,\quad \nabla^2 u=\left(
\begin{array}{cc}
u_{,xx} & u_{,xy} \\
u_{,xy} & u_{,yy} \\
\end{array}
\right)
\end{align}
and in 3D
\begin{align}
\nabla u=\Big(u_{,x},u_{,y},u_{,z}\Big)^T,\quad
\nabla^2 u=\left(
\begin{array}{ccc}
u_{,xx} & u_{,xy} & u_{,xz} \\
u_{,xy} & u_{,yy} & u_{,yz} \\
u_{,xz} & u_{,yz} & u_{,zz} \\
\end{array}
\right)
\end{align}
where $u_{,xx}$ denotes the partial derivative of $u$ with respect to $x$ twice.

In the framework of NOM, the partial derivatives can be constructed as follows. The Taylor series expansion of scalar-valued field $u_j$ in 2D can be written as
\begin{align}
u_j=u_i+(u_{i,x},u_{i,y},u_{i,xx},u_{i,xy},u_{i,yy})\cdot (x_{ij},y_{ij},x_{ij}^2/2,x_{ij} y_{ij}, y_{ij}^2/2)+O(|\br_{ij}|^3)\label{eq:ty1}
\end{align}
where $\bm r_{ij}=(x_{ij},y_{ij})^T=\bm x_j-\bm x_i$ and $O(|\br_{ij}|^3)$ denotes the higher order term.

Let
\begin{align}
u_{ij}&=u_j-u_i\\
\bm p_{ij}&=(x_{ij},y_{ij},x_{ij}^2/2,x_{ij} y_{ij}, y_{ij}^2/2)^T\\
\partial u_i&=(u_{i,x},u_{i,y},u_{i,xx},u_{i,xy},u_{i,yy})^T
\end{align}

The Taylor series expansion of Eq.\ref{eq:ty1} can be rewritten as
\begin{align}
u_{ij}=\partial u_i^T\bm p_{ij}\label{eq:ty2}
\end{align}
Tensor product with $\bm p_{ij}^T$ on both sides of Eq.\ref{eq:ty2}
\begin{align}
u_{ij}\bm p_{ij}^T =\partial u_i^T \bm p_{ij}\bm p_{ij}^T
\end{align}
Considering the weighted integration in the support $\cS_i$, we obtain
\begin{align}
\int_{\cS_i}\omega(\br_{ij}) u_{ij} \bm p_{ij}^T \ud V_j=\partial u_i^T \int_{\cS_i} \omega(\br_{ij})\bm p_{ij}\bm p_{ij}^T \ud V_j\label{eq:nomeq}
\end{align}
where $\omega(\br_{ij})$ is the weight function.

Then the nonlocal operators can be obtained as
\begin{align}
\tpartial u_i:= \int_{\cS_i}\omega(\br_{ij})\bm K_i\cdot \bm p_{ij}u_{ij} \ud V_j \label{eq:pui2}
\end{align}
where
\begin{align}
\bm K_i=\Big(\int_{\cS_i} \omega(\br_{ij})\bm p_{ij}\otimes \bm p_{ij}^T \ud V_j\Big)^{-1}\label{eq:shapeT}
\end{align}
Here, we use $\tilde{\square}$ to denote the nonlocal form of the local operator $\square$ since the definitions of the local operator and the nonlocal operator are distinct.

The Taylor series expansion of a vector field $\bu$ can be obtained in the similar manner as
\begin{align}
\bm u_{ij}^T&= \bm p_{ij}^T\cdot \partial \bm u_i\\
\omega(\br_{ij}) \bm p_{ij}\otimes \bm u_{ij}^T &=\omega(\br_{ij})\bm p_{ij}\otimes \bm p_{ij}^T \cdot \partial \bm u_i \\
\int_{\cS_i} \omega(\br_{ij}) \bm p_{ij}\otimes \bm u_{ij}^T\ud V_j &=\int_{\cS_i} \omega(\br_{ij})\bm p_{ij}\otimes \bm p_{ij}^T \cdot \partial \bm u_i\ud V_j
\end{align}
That is
\begin{align}
\tpartial \bm u_i:=\int_{\cS_i}\omega(\br_{ij})\bm K_i\cdot \bm p_{ij} \otimes \bm u_{ij}^T\ud V_j
\end{align}

For example, consider the displacement field $\bm u=(u,v)^T$ in two dimensional space, the relative displacement vector and the nonlocal partial derivatives have the explicit forms
\begin{align}
\bm u_{ij}=
\begin{pmatrix}
u_{j}-u_{i}\\
v_{j}-v_{i}\\
\end{pmatrix},
\quad
\tpartial \bm u_{i}=(\tpartial u_i, \tpartial v_i)=
\begin{pmatrix}
u_{i,x} & v_{i,x}\\
u_{i,y} & v_{i,y}\\
u_{i,xx} & v_{i,xx}\\
u_{i,xy} & v_{i,xy}\\
u_{i,yy} & v_{i,yy}\\
\end{pmatrix},
\end{align}
%and the shape tensor
%polynomial vector based on relative coordinate $\bm r_{ij}=(x_{ij},y_{ij})^T=\bm x_j-\bm x_i$ be written as
%\begin{align}
%\bm p_{ij}=(x_{ij},y_{ij},x_{ij}^2/2,x_{ij} y_{ij}, y_{ij}^2/2)^T
%\end{align}

%Therefore
%\begin{align}
%\partial \bm u_i =\int_{\cS_i}\omega(\br_{ij})\bm K_i \cdot \bm p_{ij} \otimes \bm u_{ij}^T \ud V_j
%\end{align}
Let $\bm K_i \cdot \bm p_{ij}$ be denoted by
\begin{align}
(g_{1j},g_{2j},h_{1j},h_{2j},h_{3j})^T= \bm K_i \cdot \bm p_{ij}\label{eq:p2d}
\end{align}
The gradient vector $\bm g_{ij}$ and Hessian matrix $\bm h_{ij}$ between points $i$ and $j$ in 2D are, respectively
\begin{align}
\bm g_{ij}=(g_{1j},g_{2j})^T,\quad
\bm h_{ij}=\begin{pmatrix}
h_{1j} & h_{2j}\\
h_{2j} & h_{3j}\\
\end{pmatrix}
\end{align}

In 3D case, the polynomial vector based on relative coordinates $\bm r_{ij}=(x_{ij},y_{ij},z_{ij})^T=\bm x_j-\bm x_i$ is given as
\begin{align}
\bm p_{ij}=(x_{ij},y_{ij},z_{ij},x_{ij}^2/2,x_{ij} y_{ij},x_{ij}z_{ij}, y_{ij}^2/2,y_{ij}z_{ij},z_{ij}^2)^T\label{eq:p3d}
\end{align}
The shape tensor in 3D is constructed by Eq.\ref{eq:shapeT} with $\bm p_{ij}$ in Eq.\ref{eq:p3d}.
%\begin{align}
%\bm K_i=\Big(\int_{\cS_i} \omega(\br_{ij})\bm p_{ij}\otimes \bm p_{ij}^T \ud V_j\Big)^{-1}
%\end{align}

%Therefore
%\begin{align}
%\partial \bm u_i =\int_{\cS_i}\omega(\br_{ij})\bm K_i \cdot \bm p_{ij} \otimes \bm u_{ij}^T \ud V_j
%\end{align}
Let $\bm K_i \cdot \bm p_{ij}$ in 3D be denoted by
\begin{align}
(g_{1j},g_{2j},g_{3j},h_{1j},h_{2j},h_{3j},h_{4j},h_{5j},h_{6j})^T= \bm K_i \cdot \bm p_{ij}
\end{align}

The gradient vector $\bm g_{ij}$ and Hessian matrix $\bm h_{ij}$ for two points $i,j$ in support in 3D are, respectively
\begin{align}
\bm g_{ij}=(g_{1j},g_{2j},g_{3j})^T,\quad
\bm h_{ij}=\begin{pmatrix}
h_{1j} & h_{2j} & h_{3j}\\
h_{2j} & h_{4j} & h_{5j}\\
h_{3j} & h_{5j} & h_{6j}\\
\end{pmatrix}
\end{align}
It is worth mentioning that for first order NOM or peridynamics, the gradient vector can be calculated as well by
\begin{align}
\bm g_{ij}=\Big(\int_{\cS_i} \omega(\br_{ik}) \br_{ik}\otimes \br_{ik} \ud V_k\Big)^{-1}\cdot \br_{ij}
\end{align}

Then the nonlocal gradient operator and Hessian operator for vector field can be defined as
\begin{align}
\tnabla \otimes \bm u_i&:=\int_{\cS_i} \omega(\br_{ij})\bm u_{ij} \otimes \bm g_{ij} \ud V_j\label{eq:nu1}\\
\tnabla\otimes \tnabla \otimes \bm u_i&:=\int_{\cS_i} \omega(\br_{ij})\bm u_{ij} \otimes \bm h_{ij} \ud V_j
\end{align}

In the case of 2-vector in 2 dimensional space, the explicit forms of $\tnabla \otimes \bm u_i$ and $\tnabla \otimes\tnabla \otimes \bm u_i$ are
\begin{align}
\tnabla \otimes \bm u_i=\begin{pmatrix}
u_{i,x} & u_{i,y}\\
v_{i,x} & v_{i,y}\\
\end{pmatrix}
\end{align}
\begin{align}
\tnabla\otimes \tnabla \otimes \bm u_i=
\begin{pmatrix}
\pfrac{(\tnabla \otimes \bm u_i)}{x} &
\pfrac{(\tnabla \otimes \bm u_i)}{y}\\
\end{pmatrix}=
\begin{pmatrix}
\begin{pmatrix}
u_{i,xx} & u_{i,yx}\\
v_{i,xx} & v_{i,yx}\\
\end{pmatrix} &
\begin{pmatrix}
u_{i,xy} & u_{i,yy}\\
v_{i,xy} & v_{i,yy}\\
\end{pmatrix}
\end{pmatrix}
\end{align}

For scalar-valued field, the nonlocal Laplace operator is the tensor contraction of $\tnabla \otimes \tnabla u_i$, e.g. $\tilde{\Delta}=\tnabla\cdot \tnabla=\mbox{tr}(\tnabla \otimes \tnabla)$, where $tr(\cdot)$ denotes the trace of a matrix. More specifically, in 2D
\begin{align}
\tilde{\Delta} u_i:=\int_{\cS_i} \omega(\br_{ij})(h_{1j}+2 h_{2j}+h_{3j}) u_{ij}\ud V_j\label{eq:lap2D}
\end{align}
and in 3D
\begin{align}
\tilde{\Delta} u_i:=\int_{\cS_i} \omega(\br_{ij})(h_{1j}+2 h_{2j}+2 h_{3j}+h_{4j}+2 h_{5j}+h_{6j}) u_{ij}\ud V_j\label{eq:lap3D}
\end{align}
And their local counterparts for scalar-valued field are

\begin{align}
\Delta w&= w_{,yy}+2 w_{,xy}+w_{,xx} &\mbox{ in 2D} \\
\Delta w&= w_{,xx}+w_{,yy}+w_{,zz}+2 w_{,xy}+2 w_{,xz}+2 w_{,yz}&\mbox{ in 3D}
\end{align}

\subsection{Stability of the second-order nonlocal operators}
%\begin{subequations}
%\begin{align}
%\bm K_i&=\Big(\sum_{j\in\cS_i}\omega(\br)\bp_{j}\otimes(\bp_{j})^T \Delta V_{j}\Big)^{-1},\label{eq:ibK}\\
%\bp^h_{wi}&=\Big(w(\br_{j_1}) \bp_{j_1}\Delta V_{j_1},...,w(\br_{j_{n_i}}) \bp_{j_{n_i}}\Delta V_{j_{n_i}}\Big)\label{eq:bpw}\\
%\Delta \bu_i&=(u_{i j_1},...,u_{i j_k},...,u_{i j_{n_i}})^T\label{eq:dbu}
%\end{align}
%\end{subequations}

According to Ref \cite{ren2020higher}, the energy functional for second-order nonlocal operator in discrete form can be written as
\begin{align}
\fF_i(\bu)=\frac 12 \frac{p^{hg}}{m_i} \int_{\cS_i} \omega(\br_{ij})\big(u_{ij}-\bp_j^T\tpartial u_i \big)^2 \ud V_j\label{eq:Bu}
\end{align}
where $p^{hg}$ is the penalty and $m_i=\int_{\cS_i} \omega(\br_{ij})\ud V_j$.
The operator in Eq.\ref{eq:pui2} corresponds to the minimum of Eq.\ref{eq:Bu}. The first variation of $\fF_i$ is
\begin{align}
\delta\fF_i(\bu)&=\frac{p^{hg}}{m_i} \int_{\cS_i} \omega(\br_{ij})\big(u_{ij}-\bp_j^T\tpartial u_i \big) (\delta u_j-\delta u_i-\bm p_j^T \tpartial \delta u_i) \ud V_j\notag\\
&=\frac{p^{hg}}{m_i} \int_{\cS_i} \omega(\br_{ij})\big(u_{ij}-\bp_j^T\tpartial u_i \big) (\delta u_j-\delta u_i) \ud V_j\notag\\
&-\frac{p^{hg}}{m_i}\int_{\cS_i} \omega(\br_{ij})\big(u_{ij}-\bp_j^T\tpartial u_i \big) (\bm p_j^T \tpartial \delta u_i) \ud V_j
\end{align}
We can prove that
\begin{align}
&-\frac{p^{hg}}{m_i}\int_{\cS_i} \omega(\br_{ij})\big(u_{ij}-\bp_j^T\tpartial u_i \big) (\bm p_j^T \tpartial \delta u_i) \ud V_j\notag\\
=&-\frac{p^{hg}}{m_i}\int_{j \in\cS_i} \omega(\br_{ij})\big(\bm p_j u_{ij} -\bp_j \bm p_j^T\tpartial u_i \big) \ud V_j \cdot \tpartial \delta u_i \notag\\
%=-\Big(\sum_{j \in\cS_i} \omega(\br)\bm p_j u_{ij}\Delta V_j -\sum_{j \in\cS_i} \omega(\br) \bp_j \bm p_j^T\tpartial u_i \Delta V_j\Big) \cdot \tpartial \delta u_i \Delta V_i\notag\\
=&-\frac{p^{hg}}{m_i}\Big(\underbrace{\int_{\cS_i} \omega(\br_{ij})\bm p_j u_{ij}\ud V_j -\int_{\cS_i} \omega(\br) \bp_j \bm p_j^T \ud V_j\cdot \tpartial u_i}_{=0 \text{ since Eq.\ref{eq:nomeq}}}\Big) \cdot \tpartial \delta u_i \notag\\
=&0\notag
\end{align}

Therefore,
\begin{align}
\delta\fF_i(\bu)=\frac{p^{hg}}{m_i} \int_{\cS_i} \omega(\br_{ij})\big(u_{ij}-\bp_j^T\tpartial u_i \big) (\delta u_j-\delta u_i) \ud V_j\notag
\end{align}
Consider integration of $\delta\fF_i(\bu)$ in domain
\begin{align}
&\int_{\Omega} \delta \fF_i \ud V_i=p^{hg}\underbrace{ \int_{\Omega} \int_{\cS_i}\frac{\omega(\br_{ij})}{m_i} \big(u_{ij}-\bp_j^T\tpartial u_i \big) (\delta u_j-\delta u_i) \ud V_j\ud V_i}_{\text{by Eq.\ref{eq:dualprop}}}\notag\\
=&\int_{\Omega} \Big(\int_{\cS_i'} \omega(\br_{ij})\frac{p^{hg}}{m_j}\big(u_{ji}-\bp_i^T\tpartial u_j \big)\ud V_j-\int_{\cS_i} \omega(\br)\frac{p^{hg}}{m_i}\big(u_{ij}-\bp_j^T\tpartial u_i\big) \ud V_j\Big)\delta u_i\ud V_i
\end{align}
For any $\delta u_i$, $\int_{\Omega} \delta \fF_i \ud V_i=0$ leads to the internal force due to the stability of the nonlocal operator
\begin{align}
\int_{\cS_i'} \omega(\br_{ij})\frac{p^{hg}}{m_j}\big(u_{ji}-\bp_i^T\tpartial u_j \big)\ud V_j-\int_{\cS_i} \omega(\br)\frac{p^{hg}}{m_i}\big(u_{ij}-\bp_j^T\tpartial u_i\big) \ud V_j \label{eq:nomhgij}
\end{align}
Eq.\ref{eq:nomhgij} is the expression for a scalar-valued field. For vector-valued field, the internal force due to the stability of nonlocal operator is
\begin{align}
\int_{\cS_i'} \omega(\br_{ji})\frac{p^{hg}}{m_j}\big(\bm u_{ji}-\bp_i^T\tpartial \bm u_j \big)\ud V_j-\int_{\cS_i} \omega(\br_{ij})\frac{p^{hg}}{m_i}\big(\bm u_{ij}-\bp_j^T\tpartial \bm u_i\big) \ud V_j \label{eq:nomhgijv}
\end{align}

\section{Nonlocal governing equations based on NOM}\label{sec:EqbyNOM}
This section is devoted to the variational derivation of nonlocal strong forms of solid mechanics, including hyperelasticity, thin plate, gradient elasticity, electro-magnetic-elasticity theory and phase field fracture method. The strong form is suitable for theoretical analysis as well as explicit time integration. For the fully implicit simulation of various PDEs, the reader is referred to NOM for PDEs \cite{ren2020nonlocal,rabczuk2019nonlocal,ren2020higher,ren2020nonlocalni,ren2020nonlocalCH}.

\subsection{Nonlocal form for hyperelasticity}
Consider the energy density of a hyperelasticity as $\phi:=\phi(\bm F)$, where $\bm F=\nabla \bm u+\bm I$.
The balance equation for the hyperelastic solid is
\begin{align}
\nabla\cdot \bm P+\bm b=0 \mbox{ on } \Omega \label{eq:hpe}
\end{align}
with boundary conditions $\bm u=\bm u_0 \mbox{ on } \Gamma_D$ and $\bm P\cdot \bm n=\bm t_0 \mbox{ on }\Gamma_N$, where $\bm u_0$ is the specified displacement and $\bm t_0$ is the prescribed traction load, $\bm P=\pfrac{\phi}{\bm F}$, the first Piola-Kirchhoff stress, $\bm b$ is the body force density.
\subsubsection{Derivation based on variational principle}
The variation of strain energy over the domain is
\begin{align}
\delta \cF&=\int_{\Omega} \delta \phi(\bm F) \ud V=\int_{\Omega} \pfrac{\phi}{\bm F} :\delta \bm F\ud V\notag\\
&=\int_{\Omega} \bm P :\nabla (\delta \bm u)\ud V\notag\\
&=\int_{\Omega} \bm P_i :\int_{\cS_i} \omega(\br_{ij}) \delta\bm u_{ij}\otimes \bm g_{ij} \ud V_j\ud V_i\notag\\
&=\int_{\Omega} \int_{\cS_i} \omega(\br_{ij}) \bm P_i :\delta\bm u_{ij}\otimes \bm g_{ij} \ud V_j\ud V_i\notag\\
&=\int_{\Omega} \int_{\cS_i} \omega(\br_{ij}) ( \bm P_i\cdot \bm g_{ij})\cdot \delta\bm u_{ij} \ud V_j\ud V_i\notag\\
&=\underbrace{\int_{\Omega} \int_{\cS_i} \omega(\br_{ij}) ( \bm P_i\cdot \bm g_{ij}) \cdot (\delta\bm u_{j}-\delta \bm u_i) \ud V_j\ud V_i}_{\text{by Eq.\ref{eq:dualprop}}}\notag\\
&=\int_{\Omega} \Big(\int_{\cS_i'} \omega(\br_{ji})\bm P_j \cdot \bm g_{ji} \ud V_j-\int_{\cS_i} \omega(\br_{ij}) \bm P_i \cdot \bm g_{ij} \ud V_j\Big)\cdot\delta \bu_i\ud V_i
\end{align}
In above derivation, we replace the gradient operator with nonlocal gradient, e.g. $\tnabla \otimes \bm u_i\to\int_{\cS_i} \omega(\br_{ij})\bm u_{ij} \otimes \bm g_{ij} \ud V_j$ in Eq.\ref{eq:nu1}, and the relation $\bm A:\bm a\otimes \bm b= (\bm A\cdot\bm b)\cdot\bm a$ for second-order tensor $\bm A$ and vectors $\bm a,\bm b$ is employed.

The variational of external body force energy
\begin{align}
\delta \cF_{ext}=\int_{\Omega} \delta \bm u \cdot \bm b \ud V
\end{align}
For any $\delta\bm u_i$, $\delta \cF-\delta \cF_{ext}=0$ leads to the nonlocal governing equations for elasticity %expression for $\nabla\cdot \bm P$
\begin{align}
\int_{\cS_i} \omega(\br_{ij}) \bm P_i \cdot \bm g_{ij}\ud V_j-\int_{\cS_i'} \omega(\br_{ji}) \bm P_j \cdot \bm g_{ji}\ud V_j +\bm b=0 \label{eq:dsnome}
\end{align}
%This is the governing equation of nonlocal elasticity, derived from the variational principle based on the nonlocal operator method.

Considering the effect of inertial force $\rho \ddot{\bm u}_i$ per unit volume, and replacing the dual-support with dual-horizon, we obtain the equations of motion for dual-horizon peridynamics
\begin{align}
\int_{\mathcal H_i} \omega(\br_{ij}) \bm P_i \cdot \bm g_{ij}\ud V_j-\int_{\mathcal H_i'} \omega(\br_{ji}) \bm P_j \cdot \bm g_{ji}\ud V_j +\bm b_i=\rho \ddot{\bm u}_i \label{eq:dhpd}
\end{align}
If the sizes of horizons for all material points are the same, the dual-horizon peridynamics degenerates to the conventional constant horizon peridynamics.

For any specific strain energy density (for example, isotropic/anisotropic linear/nonlinear elasticity), the explicit form of $\bm P$ can be derived straightforwardly.

\subsubsection{Derivation based on weighted residual method}
Beside the derivation based on strain energy density, the nonlocal strong form can be derived by weighted residual method.
Consider the governing equations for hyperelasticity
, the weak form of Eq.\ref{eq:hpe} for any trial vector becomes
\begin{align}
0&=\int_{\Omega} \bm v\cdot \nabla \cdot \bm P+\bm v\cdot \bm b \ud V\notag\\
&=\int_{\Omega} -\nabla \bm v: \bm P+\bm v\cdot \bm b\ud V+\int_{\Gamma} \bm P \cdot \bm n \cdot \bm v\ud S\notag\\
&=\int_{\Omega} -\Big(\int_{\cS_i} \omega(\br_{ij})\bm v_{ij} \otimes \bm g_{ij} \ud V_j\Big): \bm P_i+\bm v_i\cdot \bm b \ud V_i+\int_{\Gamma} \bm P \cdot \bm n \cdot \bm v \ud S
\end{align}
Let us focus on the integral in $\Omega$, the first term in above equation can be written as
\begin{align}
&\int_{\Omega} -\Big(\int_{\cS_i} \omega(\br_{ij})\bm v_{ij} \otimes \bm g_{ij} \ud V_j\Big): \bm P_i \ud V_i\notag\\
&= \underbrace{\int_{\Omega} -\Big(\int_{\cS_i} \omega(\br_{ij}) \bm P_i\cdot \bm g_{ij} \cdot (\bm v_{j}-\bm v_i) \ud V_j\Big) \ud V_i}_{\text{by Eq.\ref{eq:dualprop}}}\notag\\
&=\int_{\Omega} \Big(\int_{\cS_i} \omega(\br_{ij}) \bm P_i\cdot \bm g_{ij} \ud V_j-\int_{\cS_i'} \omega(\br_{ji}) \bm P_j\cdot \bm g_{ji}\ud V_j\Big)\cdot \bm v_i \ud V_i
\end{align}
For any $\bm v_i$, the weak form being zero leads to
\begin{align}
\int_{\cS_i} \omega(\br_{ij}) \bm P_i\cdot \bm g_{ij} \ud V_j-\int_{\cS_i'} \omega(\br_{ji}) \bm P_j\cdot \bm g_{ji}\ud V_j+\bm b=0\notag
\end{align}
which is identical to Eq.\ref{eq:dsnome}. As being more general than the energy method, the weighted residual method can be used to convert PDEs that have no energy functional to nonlocal integral forms.

%When $\bm g_{ij}$ is derived from the higher order NOM, these formulation can be viewed as

\subsection{Nonlocal thin plate theory}
The thin plate theory is widely used in engineering applications \cite{timoshenko1959theory}. The basic assumption of thin plate include: 1) the thickness of the plate is much smaller than the length inside the mid-plane; 2) the deflection is much smaller than the thickness of the plate so that higher order effect is neglect-able; 3) the stress along the thickness direction is assumed as zero, e.g. $\sigma_z\approx 0$ and the points in the midplane have no displacement parallel to the midplane, e.g. $u(x,y,0)=v(x,y,0)\approx 0$; 4) the normal of the mid-plane remains perpendicular to the mid-plane after deformation. Then the plate bending can be simplified into 2D problem and the displacements, strain and stress can be described by the deflection on the mid-plane
\begin{align}
u(x, y, z)&=-z \frac{\partial w}{\partial x} \\
v(x, y, z)&=-z \frac{\partial w}{\partial y} \\
w(x, y, z) &\simeq w(x, y, 0) \cong w(x, y)
\end{align}
The generalized strain is the Hessian operator on the deflection

\begin{align}
\bm \kappa=\nabla^2 w=\begin{pmatrix}w_{,xx}& w_{,xy}\\ w_{,xy}&w_{,yy}\end{pmatrix}
\end{align}
with nonlocal correspondence and its variation
\begin{align}
\bm \kappa=\tnabla^2 w:=\int_{\cS_i} \omega(\br_{ij}) \bm h_{ij} w_{ij} \ud V_j
\end{align}

\begin{align}
\delta \bm \kappa=\int_{\cS_i} \omega(\br_{ij}) \bm h_{ij} \delta w_{ij} \ud V_j
\end{align}
The momentum tensor $\bm M$, the general stress for isotropic thin plate, is given by
\begin{align}
\bm{M}=\left(\begin{array}{cc}
M_{xx} & M_{xy} \\
M_{xy} & M_{yy}
\end{array}\right)=D_{0}\big(\nu \mbox{ tr} (\bm{\kappa}) \bm{I}_{2 \times 2}+(1-\nu) \bm{\kappa}\big)
\end{align}
where $D_{0}=\frac{E t^{3}}{12\left(1-\nu^{2}\right)}$ and $t$ is the thickness of the plate.

Based on the principle of minimum potential energy, the energy functional for the governing equation is
\begin{align}
\mathcal F_{int}=\int_{\Omega} \frac 12 \bm M:\bm \kappa -q w \ud S
\end{align}
and for the boundary condition can be expressed as
\begin{align}
\mathcal F_{ext}=\int_{S_3} \bar{V}_n w\ud \Gamma-\int_{S_2+S_3} \bar{M}_n \pfrac{w}{n} \ud \Gamma
\end{align}
where $q$ is the external transverse load on the mid-plane, $\bar{V}_n$ is the shear force load on boundary $S_3$ and $\bar{M}_n$ is the prescribed moment on boundary $S_2+S_3$. For simplicity, we leave the integral on the boundary for later consideration. The variation of the internal energy functional is
\begin{align}
\delta \mathcal F_{int}&=\int_{\Omega} \bm M: \delta \bm \kappa-q \delta w \ud S\notag\\
&=\int_{\Omega} \bm M_i: \int_{\cS_i} \omega(\br_{ij}) \bm h_{ij} \delta w_{ij} \ud S_j -q_i \delta w_i\ud S_i\notag\\
&=\underbrace{\int_{\Omega}\int_{\cS_i} \omega(\br_{ij}) \bm M_i: \bm h_{ij} (\delta w_{j}-\delta w_i) \ud S_j}_{\text{by Eq.\ref{eq:dualprop}}} -\int_{\Omega}q_i \delta w_i \ud S_i\notag\\
&=\int_{\Omega}\Big( \int_{\cS_i'} \omega(\br_{ij}) \bm M_j: \bm h_{ji}\ud S_j- \int_{\cS_i} \omega(\br_{ij}) \bm M_i: \bm h_{ij} \ud S_j-q_i\Big) \delta w_i \ud S_i
\end{align}
The variation of the external energy function is
\begin{align}
\delta\mathcal F_{ext}&=\int_{S_3} \bar{V}_n \delta w\ud \Gamma-\int_{S_2+S_3} \bar{M}_n \pfrac{\delta w}{n} \ud \Gamma\notag\\
&=\int_{S_3} \bar{V}_n \delta w\ud \Gamma-\int_{S_2+S_3} \bar{M}_n \nabla \delta w\cdot \bm n \ud \Gamma\notag\\
&=\int_{S_3} \bar{V}_n \delta w\ud \Gamma-\int_{S_2+S_3} \bar{M}_{ni} \int_{\cS_i}\omega(\br_{ij}) \delta w_{ij} \bm g_{ij} \ud V_j \cdot \bm n_i \ud \Gamma_i\notag\\
&=\int_{S_3} \bar{V}_n \delta w\ud \Gamma-\int_{S_2+S_3} \int_{\cS_i}\omega(\br_{ij})\bar{M}_{ni} \bm g_{ij} \cdot \bm n_i \delta w_{ij}\ud V_j \ud \Gamma_i\notag\\
&=\int_{S_3} \bar{V}_n \delta w\ud \Gamma-\int_{S_2+S_3} \Big( \int_{\cS_i'}\omega(\br_{ji})\bar{M}_{nj} \bm g_{ji} \cdot \bm n_j \ud V_j -\int_{\cS_i}\omega(\br_{ij})\bar{M}_{ni} \bm g_{ij} \cdot \bm n_i \ud V_j\Big)\delta w_i \ud \Gamma_i
\end{align}

For any $\delta w_i$, $\delta \mathcal F_{int}-\delta \mathcal F_{ext}=0$ leads to the nonlocal thin plate equation for material point in domain $\Omega$
\begin{align}
\int_{\cS_i} \omega(\br_{ij}) \bm M_i: \bm h_{ij} \ud V_j- \int_{\cS_i'} \omega(\br_{ij}) \bm M_j: \bm h_{ji}\ud V_j+q_i=0 \label{eq:nlthinplate}
\end{align}
The additional nonlocal form for material point applied with the moment boundary condition is
\begin{align}
\int_{\cS_i}\omega(\br_{ij})\bar{M}_{ni} \bm g_{ij} \cdot \bm n_i \ud V_j- \int_{\cS_i'}\omega(\br_{ji})\bar{M}_{nj} \bm g_{ji} \cdot \bm n_j \ud V_j=0
\end{align}
%
%\begin{align}
%\nabla^2: \bm M_i= \int_{\cS_i'} \omega(\br_{ij}) \bm M_j: \bm h_{ji}\ud V_j- \int_{\cS_i} \omega(\br_{ij}) \bm M_i: \bm h_{ij} \ud V_j
%\end{align}

Based on the D'Alembert's principle, the equation of motion considering the effect of inertial force $\rho t\ddot{w}_i$ per unit area is
\begin{align}
%\tnabla^2: \bm M_i=
\int_{\cS_i'} \omega(\br_{ij}) \bm M_j: \bm h_{ji}\ud V_j- \int_{\cS_i} \omega(\br_{ij}) \bm M_i: \bm h_{ij} \ud V_j+q_i=t\rho\ddot{w}_i\label{eq:dstp}
\end{align}
For clamped boundary condition $w_{,n}=\nabla w\cdot \bm n=0$, the nonlocal form is
\begin{align}
\int_{\cS_i} \omega(\br_{ij}) w_{ij} \bm g_{ij}\cdot \bm n_i\ud V_j=0
\end{align}

Compared with the local governing equation for thin plate $\nabla^2:\bm M+q=t\rho \ddot{w}$, we can find the correspondence between local and nonlocal formulation
\begin{align}
\nabla^2:\bm M\to \tnabla^2: \bm M_i:= \int_{\cS_i'} \omega(\br_{ji}) \bm M_j: \bm h_{ji}\ud V_j- \int_{\cS_i} \omega(\br_{ij}) \bm M_i: \bm h_{ij} \ud V_j
\end{align}

The nonlocal derivation for thin plate can be extended to composite plate and functional gradient plate theories.
\subsection{Nonlocal gradient elasticity}
%\subsubsection{Review of gradient elasticity}
%Gradient elasticity theory introduces an internal length scale and higher-order gradients of the displacement field to account for size effects at the micro- or nano-scale.
Gradient theories emerge from considerations of the microstructure in the material at micro-scale, where a mass point after homogenization is not the center of a micro-volume and the rotation of the micro-volume depends on the moment stress/couple stress as well as the Cauchy stress.
%physical background of couple stress. History and development, physical meanings.
%Different from classical elasticity theory, such consideration enables gradient elasticity to model some interesting phenomena (such as size effect, the stress and strain effects on surface physics, nonlocal effect at micrometer/nanometer scale). Muller and Saul \cite{muller2004elastic} reviewed the importance of surface and interface stress effects on thin films and nano-scaled structures, including the self-organization and elastic driven instabilities of nano-structures. Fischer \etal \cite{fischer2008role} studied the role of the surface energy and surface stress in phase-transforming nano-particles.
%\subsubsection{Derivation of Nonlocal gradient elasticity}
%\begin{align}
%\nabla^2\bm u:=\int_{\cS_i} \omega(\br_{ij}) \bm u_{ij} \otimes \bm h_{ij} \ud V_j
%\end{align}
Gradient elasticity generalizes the elasticity theory by employing higher order terms of the deformation gradient or the gradient of the strain tensor. Generally, the energy density functional can be assumed as $\psi:=\psi(\bm F,\nabla \bm F)=\psi(\nabla \bm u, \nabla^2 \bm u)$, where $\bm F=\nabla \bm u+\bm I$. The total potential energy in domain is
\begin{align}
\mathcal F=\int_{\Omega}\psi-\bm b\cdot \bm u \ud V %
\end{align}
The stress tensor and generalized stress tensor of first Piola-Kirchhoff type are defined as
\begin{align}
\bm \sigma= \pfrac{\psi}{\bm F}\\
\bm \Sigma= \pfrac{\psi}{\nabla \bm F}
\end{align}
The variation of the total internal energy is
\begin{align}
\delta \mathcal F&=\int_{\Omega}\Big( \pfrac{\psi}{\bm F}: \nabla\delta \bm u+ \pfrac{\psi}{\nabla\bm F}\dot{:} \nabla^2 \delta \bm u-\bm b \cdot \delta\bm u\Big) \ud V\notag\\
&=\int_{\Omega} \Big(\bm \sigma: \nabla\delta \bm u+\bm \Sigma\dot{:} \nabla^2 \delta \bm u-\bm b \cdot \delta\bm u\Big)\ud V
\end{align}
Based on the integration by parts, the local form can be derived by
\begin{align}
\delta \mathcal F
&=\int_{\partial\Omega}\Big( \bm n\cdot \bm \sigma\cdot\delta \bm u+\bm n\cdot \bm \Sigma{:} \nabla \delta \bm u \Big)\ud S-\int_{\Omega}\Big( \nabla \cdot \bm \sigma\cdot \delta \bm u+\nabla\cdot \bm \Sigma{:} \nabla \delta \bm u+\bm b\cdot \bm u \Big)\ud V\notag\\
&=\int_{\partial\Omega}\Big( \bm n\cdot \bm \sigma\cdot\delta \bm u+\bm n\cdot \bm \Sigma : \nabla \delta \bm u-\bm n \cdot \nabla \cdot \bm \Sigma \cdot \delta\bm u\Big)\ud S-\int_{\Omega} (\nabla \cdot \bm \sigma-\nabla^2 : \bm \Sigma+\bm b)\cdot \delta \bm u \ud V\label{eq:gellocal}
\end{align}
Based on D'Alembert's principle, the governing equations for dynamic gradient elasticity can be written as
\begin{align}
\nabla \cdot \bm \sigma-\nabla^2 : \bm \Sigma+\bm b=\rho \ddot{\bm u} \mbox{ in }\Omega
\end{align}
On the other hand, do the substitutions $\nabla\delta \bm u\to \int_{\cS_i} \omega(\br_{ij}) \bm g_{ij} \otimes \delta \bm u_{ij} \ud V_j, \mbox{ and } \nabla^2\delta \bm u\to \int_{\cS_i} \omega(\br_{ij}) \bm h_{ij} \otimes \delta \bm u_{ij} \ud V_j$, we get
\begin{align}
\delta \mathcal F
&=\int_{\Omega} \bm \sigma: \nabla\delta \bm u+ \bm \Sigma\dot{:} \nabla^2 \delta \bm u -\bm b\cdot\delta \bm u \ud V\notag \\
&=\int_{\Omega}\Big( \bm \sigma_i: \int_{\cS_i} \omega(\br_{ij}) \bm g_{ij} \otimes \delta \bm u_{ij} \ud V_j+ \bm \Sigma_i\dot{:} \int_{\cS_i} \omega(\br_{ij}) \bm h_{ij} \otimes \delta \bm u_{ij} \ud V_j-\bm b\cdot\delta \bm u\Big)\ud V_i\notag \\
&=\underbrace{\int_{\Omega} \int_{\cS_i} \omega(\br_{ij}) \bm \sigma_i {:} (\delta\bm u_{j}-\delta\bm u_{j}) \otimes\bm g_{ij} \ud V_j \ud V_i}_{\text{by Eq.\ref{eq:dualprop}}}\notag\\
&+\underbrace{\int_{\Omega} \int_{\cS_i} \omega(\br_{ij}) \bm \Sigma_i\dot{:} (\delta\bm u_{j}-\delta\bm u_{j}) \otimes\bm h_{ij} \ud V_j \ud V_i}_{\text{by Eq.\ref{eq:dualprop}}}-\int_{\Omega}\bm b\cdot\delta \bm u_i \ud V_i\notag\\
&=\int_{\Omega}\Big(\int_{\cS_i'} \omega(\br_{ji}) \bm \sigma_j\cdot\bm g_{ji}\ud V_j- \int_{\cS_i} \omega(\br_{ij}) \bm \sigma_i\cdot\bm g_{ij} \ud V_j\Big) \cdot \delta\bm u_i \ud V_i\notag\\
&+\int_{\Omega}\Big(\int_{\cS_i'} \omega(\br_{ji}) \bm \Sigma_j:\bm h_{ji}\ud V_j- \int_{\cS_i} \omega(\br_{ij}) \bm \Sigma_i:\bm h_{ij} \ud V_j\Big) \cdot \delta\bm u_i \ud V_i-\int_{\Omega}\bm b\cdot\delta \bm u_i \ud V_i\label{eq:gelnonlocal}
\end{align}

In the above derivation, we used $\bm \Sigma\dot{:}\bm u \otimes \bm h= (\bm \Sigma:\bm h) \cdot \bm u$.
%\begin{align}
%\bm \Sigma\dot{:}\bm u \otimes \bm h=\Sigma_{ijk} u_i h_{jk}= (\bm \Sigma:\bm h) \cdot \bm u
%\end{align}
For any $\delta \bm u_i$, $\delta \cF=0$ leads to the nonlocal form of gradient elasticity
\begin{align}
\int_{\cS_i} \omega(\br_{ij}) (\bm \sigma_i \cdot \bm g_{ij} +\bm \Sigma_i: \bm h_{ij})\ud V_j-\int_{\cS_i'} \omega(\br_{ji}) (\bm \sigma_j \cdot \bm g_{ji} + \bm \Sigma_j:\bm h_{ji})\ud V_j+\bm b=\rho \ddot{\bm u}_i
\end{align}
The inertia force term is added based on D'Alembert's principle.

Comparing Eq.\ref{eq:gellocal} and Eq.\ref{eq:gelnonlocal}, the correspondence from local form to nonlocal form is
\begin{align}
\nabla^2: \bm \Sigma_i \to \int_{\cS_i'} \omega(\br_{ji}) \bm \Sigma_j: \bm h_{ji}\ud V_j- \int_{\cS_i} \omega(\br_{ij}) \bm \Sigma_i: \bm h_{ij} \ud V_j
\end{align}

%\begin{align}
%L_2(\bm u)=\sqrt{\frac{\sum_{i=1}^N \Delta V_i \|\bm u_i-\bm u_i^{exact}\|^2}{\sum_{i=1}^N \Delta V_i \|\bm u_i^{exact}\|^2}}
%\end{align}
%where $ \|\bm v\|=\sqrt{\bm v\cdot \bm v}$ for vector and $ \|\bm m\|=\sqrt{\bm m: \bm m}$ for second-order tensor.

%The explicit form of nonlocal Laplace operator in 3D is
%\begin{align}
%\tilde{\nabla}^2:\bm M_i=\int_{\mathcal S_i'} \omega(\br_{ji})\bm M_j : \left(
%\begin{array}{ccc}
%h_{1i} & h_{2i} & h_{3i} \\
%h_{2i} & h_{4i} & h_{5i} \\
%h_{3i} & h_{5i} & h_{6i} \\
%\end{array}
%\right) \ud V_j-\int_{\mathcal S_i} \omega(\br_{ij})\bm M_i : \left(
%\begin{array}{ccc}
%h_{1j} & h_{2j} & h_{3j} \\
%h_{2j} & h_{4j} & h_{5j} \\
%h_{3j} & h_{5j} & h_{6j} \\
%\end{array}
%\right) \ud V_j
%\end{align}

\subsection{Nonlocal form of magneto-electro-elasticity}
In accordance with reference \cite{Liu2014Feb}, let us postulate the following form of internal energy for the energy function $\psi:=\psi(\bm F,\nabla \bm F, \bm p,\nabla\bm p, \bm m,\nabla \bm m)$, a function depends on the displacement gradient $\bm F=\nabla \bm u+\bm I$ and its second gradient $\nabla \bm F=\nabla^2 \bm u$, polarization vector $\bm p$ and its gradient $\nabla\bm p$, magnetic field $\bm m$ and its gradient $\nabla \bm m$. The total potential energy in the domain can be written as
\begin{align}
\cF=\int_{\Omega} \psi(\bm F,\nabla \bm F, \bm p,\nabla\bm p, \bm m,\nabla \bm m) \ud V
\end{align}
This model has a strong physical background, for example, the nonlinear electro-gradient elasticity for semiconductors \cite{nguyen2019nurbs} and flexoelectricity \cite{roy2019conformal}.

The first variation of $\cF$ is
\begin{align}
\delta\cF=&\int_{\Omega} \delta\psi \ud V\notag\\
=&\int_{\Omega} \pfrac{\psi}{\bm F}:\nabla \delta \bm u+\pfrac{\psi}{\nabla \bm F}\dot{:}\nabla^2 \delta \bm u+\pfrac{\psi}{\bm p}\cdot \delta\bm p+\notag\\
&\pfrac{\psi}{\nabla \bm p}:\nabla \delta \bm p+\pfrac{\psi}{\bm m}\cdot \delta\bm m+\pfrac{\psi}{\nabla \bm m}:\nabla \delta \bm m \ud V\notag\\
=&\int_{\Omega} \bm \sigma:\nabla \delta \bm u+\bm \Sigma\dot{:}\nabla^2 \delta \bm u+\bm e\cdot \delta\bm p\notag\\
&+\bm E:\nabla \delta \bm p+\bm s\cdot \delta\bm m+\bm S:\nabla \delta \bm m \ud V
\end{align}
where
\begin{align}
\bm \sigma=\pfrac{\psi}{\bm F},
\bm \Sigma=\pfrac{\psi}{\nabla \bm F},
\bm e=\pfrac{\psi}{\bm p}\\
\bm E=\pfrac{\psi}{\nabla \bm p},
\bm s=\pfrac{\psi}{\bm m},
\bm S=\pfrac{\psi}{\nabla \bm m}
\end{align}
Doing substitutions $\nabla \delta \bm u_i\to \int_{\cS_i} \omega(\br_{ij})\delta\bm u_{ij}\otimes \bm g_{ij} \ud V_j$, $\nabla^2 \delta \bm u_i\to \int_{\cS_i} \omega(\br_{ij})\delta\bm u_{ij}\otimes \bm h_{ij} \ud V_j$, $\nabla \delta \bm p_i\to \int_{\cS_i} \omega(\br_{ij})\delta\bm p_{ij}\otimes \bm g_{ij} \ud V_j$,$\nabla \delta \bm m_i\to \int_{\cS_i} \omega(\br_{ij})\delta\bm m_{ij}\otimes \bm g_{ij} \ud V_j$ and following the same operations in prior sections, the functional becomes
\begin{align}
\delta\cF&=\int_{\Omega} \Big(\int_{\cS_i'} \omega(\br_{ji})(\bm \sigma_j\cdot \bm g_{ji}+\bm \Sigma_j :\bm h_{ji}) \ud V_j-\int_{\cS_i} \omega(\br_{ij})(\bm \sigma_i\cdot \bm g_{ij}+\bm \Sigma_i :\bm h_{ij}) \ud V_j\Big)\cdot \delta \bm u_i \ud V_i\notag\\
&+\int_{\Omega} \Big(\int_{\cS_i'} \omega(\br_{ji})(\bm E_j\cdot \bm g_{ji}) \ud V_j-\int_{\cS_i} \omega(\br_{ij})\bm E_i\cdot \bm g_{ij} \ud V_j+\bm e_i\Big)\cdot \delta \bm p_i \ud V_i+\notag\\
&\int_{\Omega} \Big(\int_{\cS_i'} \omega(\br_{ji})(\bm S_j\cdot \bm g_{ji}) \ud V_j-\int_{\cS_i} \omega(\br_{ij})\bm S_i\cdot \bm g_{ij} \ud V_j+\bm s_i\Big)\cdot \delta \bm m_i \ud V_i
\end{align}
For any $\delta \bm u_i,\delta \bm p_i,\delta \bm m_i$, $\delta \cF=0$ leads to general nonlocal governing equation for mechanical field, electrical field and magnetic field, respectively
\begin{align}
\int_{\cS_i} \omega(\br_{ij})(\bm \sigma_i\cdot \bm g_{ij}+\bm \Sigma_i :\bm h_{ij}) \ud V_j-&\notag\\
\int_{\cS_i'} \omega(\br_{ji})(\bm \sigma_j\cdot \bm g_{ji}+\bm \Sigma_j :\bm h_{ji}) \ud V_j+\bm b_i&=0\\
\int_{\cS_i} \omega(\br_{ij})\bm E_i\cdot \bm g_{ij} \ud V_j-\int_{\cS_i'} \omega(\br_{ji})\bm E_j\cdot \bm g_{ji} \ud V_j-\bm e_i&=0\\
\int_{\cS_i} \omega(\br_{ij})\bm S_i\cdot \bm g_{ij} \ud V_j-\int_{\cS_i'} \omega(\br_{ji})\bm S_j\cdot \bm g_{ji} \ud V_j-\bm s_i&=0
\end{align}
In the derivation, we did not specify the exact form of the energy density, whether it is of small deformation or of finite deformation. For the specified energy form, one only needs to derive the expression for $\bm \sigma, \bm \Sigma, \bm e,\bm E,\bm s,\bm S$ based on the material constitutions.
It can be seen that the nonlocal governing equations for the continuum magneto-electro-elasticity can be obtained with ease by using nonlocal operator method and variational principle. The same rule applies for many other physical problems.

\subsection{Nonlocal form of phase field fracture method}
Phase field fracture method is powerful in fracture modelling \cite{miehe2010thermodynamically}.
The difference in tensile and compressive strengths of the material can be considered by dividing the strain energy density into a tensile part affected by the phase field and a compressive part, which is independent of the phase field,
\begin{align}
\psi_e(\bm\varepsilon(\nabla\bm u),s)=(1-s)^2\psi_e^{+} (\bm\varepsilon(\nabla\bm u))+\psi_e^{-} (\bm\varepsilon(\nabla\bm u)).
\end{align}
where $\psi_e^+ $ ($\psi_e^- $) denotes the strain energy density for tensile (compressive) part, $\bm u $ the displacement, $s\in [0,1]$ the phase field, $\bm\varepsilon $ the strain and $\ell$ is the phase field intrinsic length scale.

The full potential functional of the phase field fracture model reads
\begin{align}
\cF_\ell (\bm u,s)&=\int_\Omega\Big((1-s)^2\psi_e^{+} (\bm\varepsilon(\nabla\bm u))+\psi_e^{-} (\bm\varepsilon(\nabla\bm u))\Big)\ud V-\int_{\partial\Omega} \bm t^*\cdot\bm u\ud A\notag\\
&-\int_\Omega\bm b\cdot\bm u\ud V+\int_\Omega g_c(\frac{s^2}{2 \ell} +\frac{\ell}{2} \nabla s\cdot\nabla s)\ud V,
\end{align}
where $\bm t ^ * $ the surface traction at the boundary, $\bm b $ the body force density and $g_c$ is the critical energy release rate.

For the sake of simplicity, we neglect the surface traction force and consider the first variation of $\cF_\ell$
\begin{align}
\delta \cF_\ell&=\int_\Omega\delta \Big((1-s)^2\psi_e^{+}+\psi_e^{-}\Big)\ud V-\int_\Omega\bm b\cdot \delta\bm u\ud V+\int_\Omega g_c \delta(\frac{s^2}{2\ell} +\frac{\ell}{2} \nabla s\cdot\nabla s)\ud V\notag\\
&=\int_\Omega \Big((1-s)^2\pfrac{\psi_e^{+}}{\bm \varepsilon}:\nabla\delta \bm u-2\psi_e^{+}(1-s)\delta s +\pfrac{\psi_e^{-}}{\bm \varepsilon}:\nabla\delta \bm u\Big)\ud V-\int_\Omega\bm b\cdot \delta\bm u\ud V\notag\\&+\int_\Omega g_c (\frac{s}{\ell} \delta s+\ell\nabla s\cdot\nabla\delta s)\ud V\notag\\
&=\int_\Omega \Big(((1-s)^2\bm \sigma^++\bm \sigma^-):\nabla\delta \bm u-\bm b\cdot \delta\bm u\Big)\ud V+\int_\Omega g_c (\frac{s}{\ell} \delta s-2\frac{\psi_e^{+}}{g_c}(1-s)\delta s +\ell\nabla s\cdot\nabla\delta s)\ud V\notag\\
&=\int_\Omega \Big(\bm \sigma_i:\nabla\delta \bm u_i-\bm b_i\cdot \delta\bm u_i\Big)\ud V_i+\int_\Omega g_c (\frac{s_i}{\ell} \delta s_i-2\frac{\psi_{ei}^{+}}{g_c}(1-s_i)\delta s_i +\ell\nabla s_i\cdot\nabla\delta s_i)\ud V_i\notag\\
&=\int_\Omega \Big(\bm \sigma_i:(\int_{\cS_i} \omega(\br_{ij})\delta\bm u_{ij}\otimes \bm g_{ij} \ud V_j)-\bm b_i\cdot \delta\bm u_i\Big)\ud V_i\notag\\&+\int_\Omega g_c (\frac{s_i}{\ell} \delta s_i-2\frac{\psi_{ei}^{+}}{g_c}(1-s_i)\delta s_i +\ell\nabla s_i\cdot\int_{\cS_i} \omega(\br_{ij})\delta s_{ij} \bm g_{ij} \ud V_j)\ud V_i\notag\\
&=\int_\Omega \Big((\int_{\cS_i'} \omega(\br_{ji})\bm \sigma_j\cdot \bm g_{ji}\ud V_j-\int_{\cS_i} \omega(\br_{ij})\bm \sigma_i\cdot \bm g_{ij}\ud V_j)\cdot \delta\bm u_i-\bm b_i\cdot \delta\bm u_i\Big)\ud V_i\notag\\
&+\int_\Omega g_c \Big(\frac{s_i}{\ell} -2\frac{\psi_{ei}^{+}}{g_c}(1-s_i) +\int_{\cS_i'} \omega(\br_{ji})\ell\nabla s_j\cdot\bm g_{ji} \ud V_j-\int_{\cS_i} \omega(\br_{ij})\ell\nabla s_i\cdot\bm g_{ij} \ud V_j\Big)\delta s_i\ud V_i
\end{align}
where
\begin{align}
\bm \sigma^+=\pfrac{\psi_e^{+}}{\bm \varepsilon}, \bm \sigma^-=\pfrac{\psi_e^{-}}{\bm \varepsilon}\\
\bm \sigma=(1-s)^2\bm \sigma^++\bm \sigma^-
\end{align}
For any $\delta \bm u_i,\delta s_i$, $\delta \cF_\ell=0$ leads to the nonlocal governing equations for the mechanical field and phase field
\begin{align}
\int_{\cS_i} \omega(\br_{ij})\bm \sigma_i\cdot \bm g_{ij}\ud V_j-\int_{\cS_i'} \omega(\br_{ji})\bm \sigma_j\cdot \bm g_{ji}\ud V_j+\bm b_i&=0\\
\frac{s_i}{\ell} -2\frac{\psi_{ei}^{+}}{g_c}(1-s_i) +\int_{\cS_i'} \omega(\br_{ji})\ell\nabla s_j\cdot\bm g_{ji} \ud V_j-\int_{\cS_i} \omega(\br_{ij})\ell\nabla s_i\cdot\bm g_{ij} \ud V_j&=0
\end{align}

The above examples aim at illustrating the power of nonlocal operator method combined with weighted residual method or variational principle in the derivation of nonlocal strong forms based on their local strong or energy forms. The derived nonlocal strong forms are variationally consistent and allow variable support sizes for each point in the model. %It can be regarded as a general and efficient method to convert traditional local form to the nonlocal strong form. %This open avenue to the nonlocal integral form of porous mechanics and flexoelectricity.

\section{Instability criterion for fracture modelling}\label{sec:fracCriterion}
%Many numerical methods are developed for fracture modelling, for example, XFEM, phase field method, cohesive fracture methods, peridynamics.
Typical methods for fracture modelling are either based on diffusive crack domain in phase field methods or on direct topological modification on meshes in XFEM or bonds in PD. Direct topological modification on meshes often leads to instability issues. For example, in NOSBPD, the breakage of a bond based on the quantities derived from stress state or strain state often introduces too much perturbation to the scheme, which may abort the calculation because of the singularity in shape tensors. These criteria include critical stretch \cite{Silling2000,dipasquale2014crack}, energy based \cite{foster2011energy} or stress based criterion \cite{Zhou2016Oct,zhou2016numerical}. Another issue in NOSBPD is that the strain energy carried by a bond is not independent with other bonds. It also depends on the direction, the length of the bond, the choice of influence functions. Removing one neighbour often gives rise to catastrophic result on the calculation. A criterion on how to remove the neighbours safely from the neighbour list remains unclear.

%
%
%
%Griffith's criterion
%\begin{align}
%\sigma_f\sqrt{a} \approx \sqrt{\frac{2E\gamma}{\pi}}
%\end{align}
%
%By analogy, the critical bond strain may be written as
%\begin{align}
%\varepsilon_f=\frac{\sigma_f}{E} \approx \sqrt{\frac{\gamma}{\pi a E}}
%\end{align}
%
%This rule is problematic. The percentage of hourglass energy is much smaller than the strain energy.

%In bond-based peridynamics \cite{dipasquale2014crack}, the critical stretch between two particles depends on the critical energy release rate as
%\begin{align}
%s_{ij}^{max}=\begin{cases}\sqrt{\frac{4 \pi G_0}{9 E \delta}} \mbox{plane stress}\\
%\sqrt{\frac{5 \pi G_0}{12 E \delta}} \mbox{plane strain}\\
%\sqrt{\frac{5 G_0}{6 E \delta}} \mbox{3D}\\
%\end{cases}
%\end{align}
%where $\delta $ is the radius of the
Damage is a process deviated from the robust mathematical expression, where the transition happens in a very narrow zone, such as the crack tip front. It is observed that around the crack tip, the gradient or strain undergoes a sharp transition within a very small zone. Most conventional numerical methods for fracture modelling focus on accurate description of the singularity occurring around the crack tip, such a description is very hard to tackle and its evolution is inconvenient to update. This dilemma can be handled when something different from continuous function is introduced.

In NOM, the gradient operator is defined in a ``redundant'' way. Around the crack tip, the deformation is irregular and the part due to hourglass energy is comparable to the strain energy carried by a particle. More specifically, the operator energy in nonlocal operator method describes the irregularity of a function around the crack tip. The irregularity is the part that cannot be described by the continuous function. For continuous domain, the strain energy density is much larger than the operator energy density. However, for particles around the crack tip, the operator energy density is far from zero and the irregularity due to the singularity around the crack tip increases comparably to the strain energy density. In this sense, the operator energy density can be viewed as an indicator for the crack tip.

Unlike the strain energy density, the hourglass energy density describes the irregular deformation around the crack tip. It depends on the penalty for the strain energy. Larger penalty improves the continuity of deformation, but the extent of hourglass energy compared with the strain energy density is hard to estimate. In this paper, we propose a special manner to estimate the critical hourglass strain. Let the critical bond strain be denoted by $s_{max}$, which may depend on the characteristic length scale of the support, critical energy release rate and the elastic modulus. When the maximal strain reached $s_{max}$, the damage process is activated and the critical hourglass strain $s^{hg}_{max}$ is set as the maximal hourglass strain $s^{hg}_{ij}$ for all bonds in the computational model. In the sequential calculation, when the hourglass strain of a bond is larger than $s^{hg}_{max}$, the damage on that bond occurs, which is mathematically described as
\begin{align}
d_{ij}=\begin{cases}
0 \mbox{ if } s^{hg}_{ij}(t)>s^{hg}_{max}, t\in [0,T]\\
1 \mbox{ otherwise }
\end{cases}
\end{align}
where $d_{ij}$ denotes the damage status between particle $i$ and particle $j$.

The damage of a particle is calculated as
\begin{align}
d_{i}=\frac{\int_{\cS_i} d_{ij} \ud V_j}{\int_{\cS_i} 1 \ud V_j}
\end{align}
Every time one particle is removed from the neighbour list, the nonlocal gradient for the central particle should be recalculated based on the remaining ``healthy'' neighbour. We will apply this rule to model fractures in 2D and 3D linear elastic material. %It was found that the fracture propagates gradually without introducing much perturbation.
%
%In peridynamics, the damage happens when the stretch between two particles exceeds a certain threshold value. In analogy, we employ the same criterion in peridynamics to model fractures.
%
%We calculate the hourglass energy for each bond,
\section{Numerical implementation}
We have applied NOM to derive the nonlocal strong forms for the traditional continuum model in \S \ref{sec:EqbyNOM}. Two representive nonlocal theories, the dual-horizon peridynamics by Eq.\ref{eq:dhpd} for fracture modeling and the nonlocal thin plate by Eq.\ref{eq:dstp}, are selected for numerical test. For the DH-PD, the focus is on the test of instability criterion for quasi-static fracture modeling by explicit time integration method. The nonlocal thin plate is compared with finite element method.

The primary step in the implementation is the calculation of internal force based on the governing equations. In the first step, the computational domain is discretized into particles.
\begin{align}
\Omega=\sum_{i=1}^N \Delta V_i
\end{align}
where $N$ is the number of particles in the domain.
Then the support of each particle is represented by a list of particle indices,
\begin{align}
\cS_i=\{j_1,j_2,...,j_{n_i}\}
\end{align}
where $j$ is the global index of the particle and $n_i$ is the number of particles in $\cS_i$.

The gradient $\bm g_{ij}$ and Hessian $\bm h_{ij}$ for two particles $i,j$ can be assembled by collecting terms in $\bm K_i \cdot \bm p_{ij}$ according to Eq.\ref{eq:p2d} or Eq.\ref{eq:p3d}, where
\begin{align}
\bm K_i=\Big(\sum_{\cS_i} \omega(\br_{ij})\bm p_{ij}\otimes \bm p_{ij}^T \Delta V_j\Big)^{-1}\label{eq:shapeTd}
\end{align}
with weight function $\omega(\br_{ij})=1/|\br_{ij}|^2$.

The nonlocal differential derivatives at point $i$ can be calculated as
\begin{align}
\tpartial u_i= \sum_{j\in \cS_i}\omega(\br_{ij})\bm K_i\cdot \bm p_{ij}u_{ij} \Delta V_j \label{eq:pui}
\end{align}
The nonlocal operators in $\tpartial u_i$ can be used to define the strain tensor, stress tensor, moment and others.

In discrete form, Eq.\ref{eq:dhpd} and Eq.\ref{eq:nlthinplate} become
\begin{align}
\sum_{\mathcal H_i} \omega(\br_{ij}) \bm P_i \cdot \bm g_{ij}\Delta V_j\Delta V_i-\sum_{\mathcal H_i'} \omega(\br_{ji}) \bm P_j \cdot \bm g_{ji}\Delta V_j\Delta V_i +\bm b_i\Delta V_i=\rho \Delta V_i \ddot{\bm u}_i \label{eq:dhpdd}
\end{align}

\begin{align}
\sum_{\cS_i} \omega(\br_{ij}) \bm M_i: \bm h_{ij} \Delta V_j\Delta V_i- \sum_{\cS_i'} \omega(\br_{ij}) \bm M_j: \bm h_{ji}\Delta V_j\Delta V_i+q_i \Delta V_i=t \rho \Delta V_i \ddot{ w}_i \label{eq:nlthinplated}
\end{align}
In Eq.\ref{eq:dhpdd} and Eq.\ref{eq:nlthinplated}, the volume of particle $i$ is multiplied on both sides of the equations. It is not required to calculate the internal forces from the dual-support. Let $\bm f_i=\bm 0, 1\leq i\leq N$ denote the initial internal force on particle $i$. For each particle, one only needs to focus on the support, calculating the forces and adding the force to the particle internal force
\begin{align}
\sum_{j\in \mathcal S_i} \omega(\br_{ij}) \bm P_i \cdot \bm g_{ij}\Delta V_j\Delta V_i &\to \bm f_i
\notag\\
-\omega(\br_{ij_1}) \bm P_i \cdot \bm g_{ij_1}\Delta V_{j_1}\Delta V_i &\to \bm f_{j_1}\notag\\
-\omega(\br_{ij_2}) \bm P_i \cdot \bm g_{ij_2}\Delta V_{j_2}\Delta V_i &\to \bm f_{j_2}\notag\\
... &\notag\\
-\omega(\br_{ij_{n_i}}) \bm P_i \cdot \bm g_{ij_{n_i}}\Delta V_{j_{n_i}}\Delta V_i&\to \bm f_{j_{n_i}}
\end{align}
where $a\to b$ denotes the addition of $a$ to $b$. The process of adding force $ -\omega(\br_{ij_1}) \bm P_i \cdot \bm g_{ij}\Delta V_{j}\Delta V_i$ to $\bm f_{j}$ is equivalent to accumulating the internal forces from particle $j$'s dual-support.

For the calculating of internal force of thin plate, the same applies
\begin{align}
\sum_{j\in \mathcal S_i} \omega(\br_{ij}) \bm M_i : \bm h_{ij}\Delta V_j\Delta V_i &\to \bm f_i
\notag\\
-\omega(\br_{ij_1}) \bm M_i : \bm h_{ij_1}\Delta V_{j_1}\Delta V_i &\to \bm f_{j_1}\notag\\
-\omega(\br_{ij_2}) \bm M_i : \bm h_{ij_2}\Delta V_{j_2}\Delta V_i &\to \bm f_{j_2}\notag\\
... &\notag\\
-\omega(\br_{ij_{n_i}}) \bm M_i : \bm h_{ij_{n_i}}\Delta V_{j_{n_i}}\Delta V_i&\to \bm f_{j_{n_i}}
\end{align}

In order to maintain the stability of the nonlocal operator, the discrete form of Eq.\ref{eq:nomhgijv} is
\begin{align}
\sum_{\cS_i'} \omega(\br)\frac{p^{hg}}{m_j}\big(\bu_{ji}-\bp_i^T\tpartial \bu_j \big)\Delta V_j\Delta V_i-\sum_{\cS_i} \omega(\br)\frac{p^{hg}}{m_i}\big(\bu_{ij}-\bp_j^T\tpartial \bu_i\big) \Delta V_j \Delta V_i \label{eq:nomhgijd}
\end{align}
For particle $i$ with support $\cS_i$, the hourglass force is calculated as follows

\begin{align}
\sum_{j\in \mathcal S_i} \omega(\br_{ij})\frac{p^{hg}}{m_i}\big(\bu_{ij}-\bp_j^T\tpartial \bm u_i\big)\Delta V_{j}\Delta V_i &\to \bm f_i
\notag\\
-\omega(\br_{ij_1})\frac{p^{hg}}{m_i}\big(\bu_{ij_1}-\bp_{j_1}^T\tpartial \bm u_i\big)\Delta V_{j_{1}}\Delta V_i &\to \bm f_{j_1}\notag\\
-\omega(\br_{ij_2})\frac{p^{hg}}{m_i}\big(\bu_{ij_2}-\bp_{j_2}^T\tpartial \bm u_i\big)\Delta V_{j_{2}}\Delta V_i &\to \bm f_{j_2}\notag\\
... &\notag\\
-\omega(\br_{ij_{n_i}})\frac{p^{hg}}{m_i} \big(\bu_{ij_{n_i}}-\bp_{j_{n_i}}^T\tpartial \bm u_i\big)\Delta V_{j_{n_i}}\Delta V_i&\to \bm f_{j_{n_i}}
\end{align}

When the internal force is attained and the contribution of the external force boundary condition or body force is accumulated, the basic Verlet algorithm \cite{verlet1967computer} outlined as follows is used to update the displacement
\begin{align}
\bm u_i(t+\Delta t)=\bm u_i(t)+\bm v_i(t) \Delta t+\frac 12 \bm a_i(t)\Delta t^2\\
\bm v_i(t+\Delta t)=\bm v_i(t)+\frac 12 \Big(\bm a_i(t)+\bm a_i(t+\Delta t)\Big)\Delta t
\end{align}
where $\bm u_i$ denotes the displacement or deflection, $\bm v_i$ the velocity and $\bm a_i=\frac{\bm f_i}{m_i}$ the acceleration for particle $i$ with mass $m_i$ subject to net force $\bm f_i$. For the detailed implementation and the numerical examples, the reader can find the open source code on Github \url{https://github.com/hl-ren/Nonlocal_elasticity} and \url{https://github.com/hl-ren/Nonlocal_thin_plate}.

\section{Numerical examples}
\subsection{Accuracy of nonlocal Hessian operator}
We test two cases with analytical function
\begin{align}
w(x,y)&=x^2+y^2, \mbox{ in 2D}\\
w(x,y,z)&=x^2+y^2+z^2 \mbox{ in 3D}
\end{align}
The exact second-order partial derivatives are
\begin{align}
w_{,xx}=w_{,yy}=2, w_{,xy}=0\mbox{ in 2D}\\
w_{,xx}=w_{,yy}=w_{,zz}=2, w_{,xy}=w_{,yz}=w_{,xz}=0 \mbox{ in 3D}
\end{align}

\begin{figure}[htp]
\centering
\includegraphics[width=14cm]{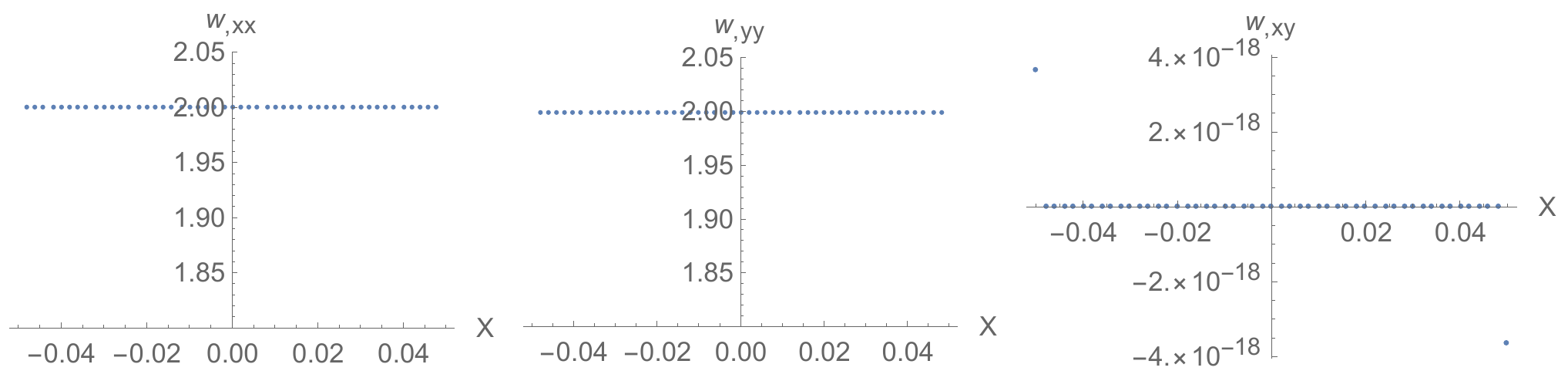}
\caption{The nonlocal Hessian with 10 neighbors in each support. }
\label{fig:N10HessianDeri}
\end{figure}

\begin{figure}[htp]
\centering
\includegraphics[width=14cm]{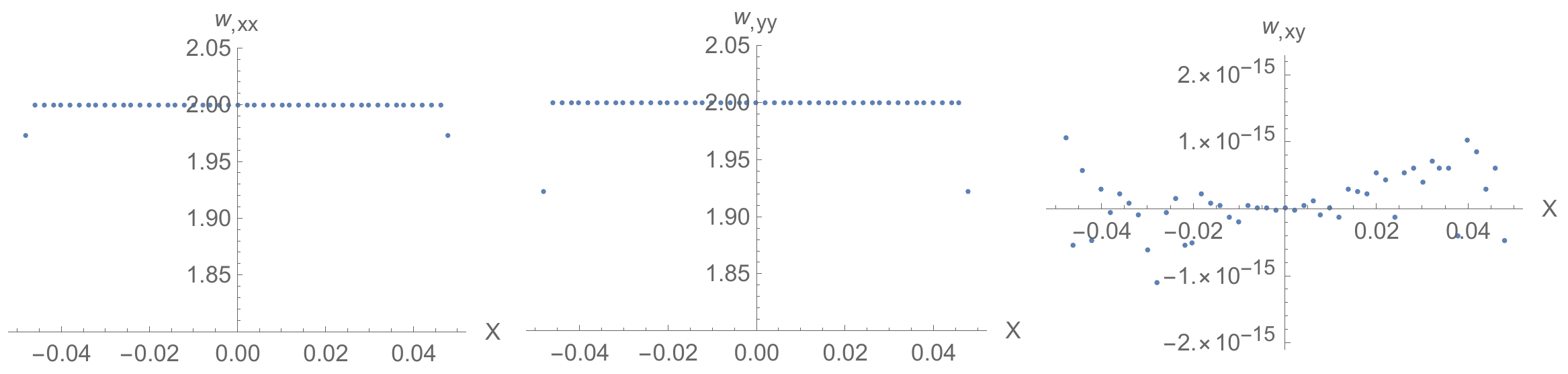}
\caption{The nonlocal Hessian with 24 neighbors in each support. }
\label{fig:N24HessianDer}
\end{figure}

\begin{figure}[htp]
\centering
\includegraphics[width=12cm]{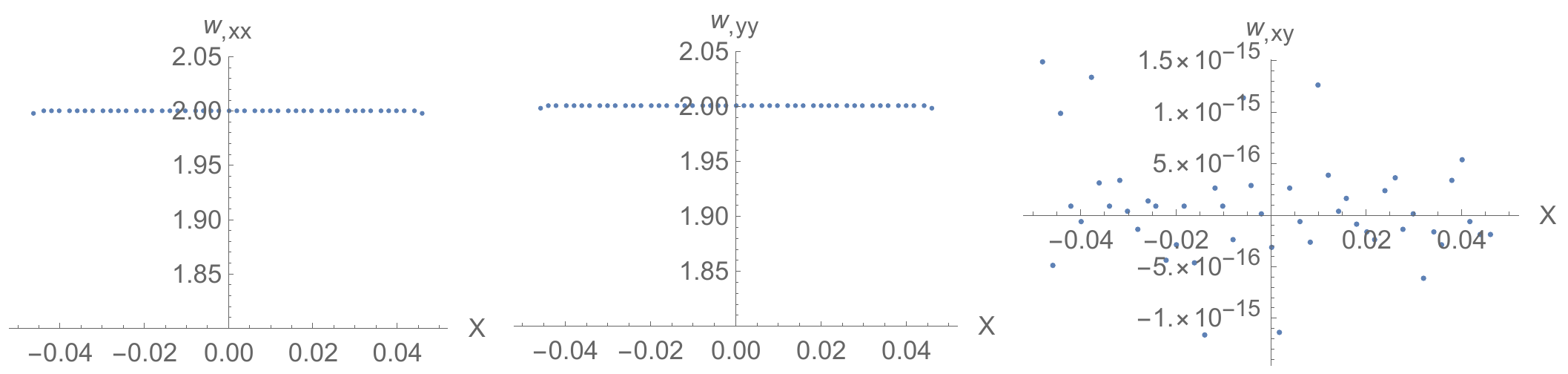}
\caption{The nonlocal Hessian with 10 neighbors in each support. }
\label{fig:N48HessianDer}
\end{figure}

\begin{figure}[htp]
\centering
\includegraphics[width=3cm]{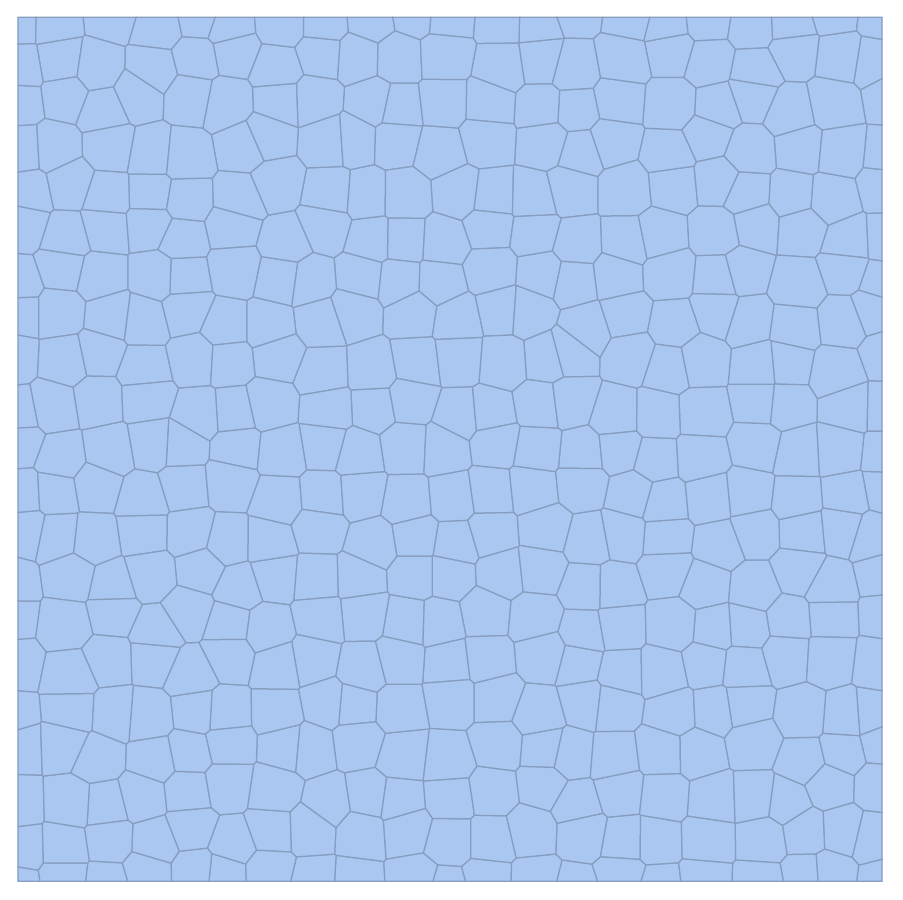}
\caption{The particle distribution based on Voronoi mesh. }
\label{fig:voronoiParticle}
\end{figure}

\begin{figure}[htp]
\centering
\includegraphics[width=12cm]{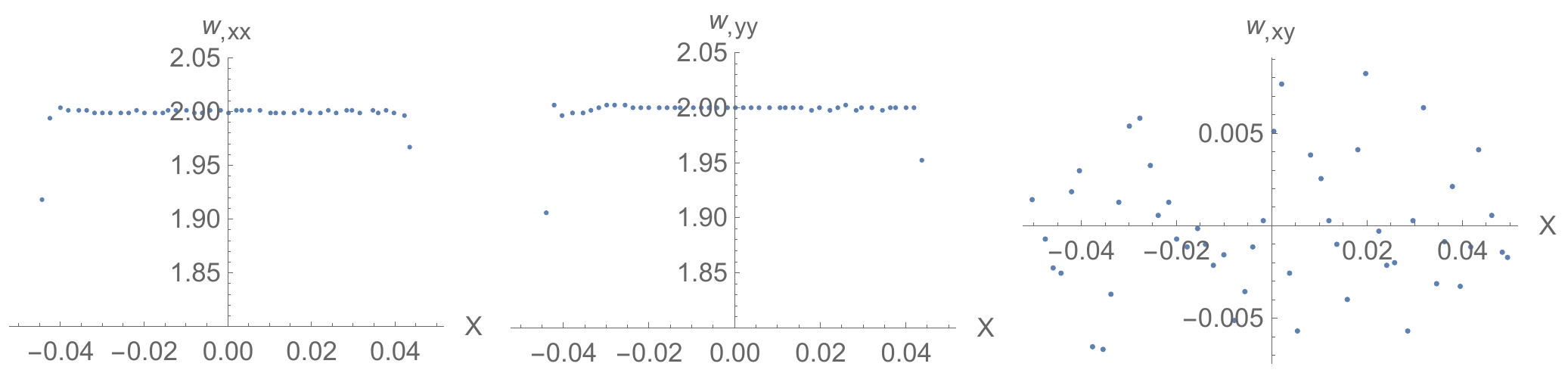}
\caption{The nonlocal Hessian with 120 neighbors in each support. }
\label{fig:N120HessianDeriInhomo}
\end{figure}

\begin{figure}[htp]
\centering
\includegraphics[width=12cm]{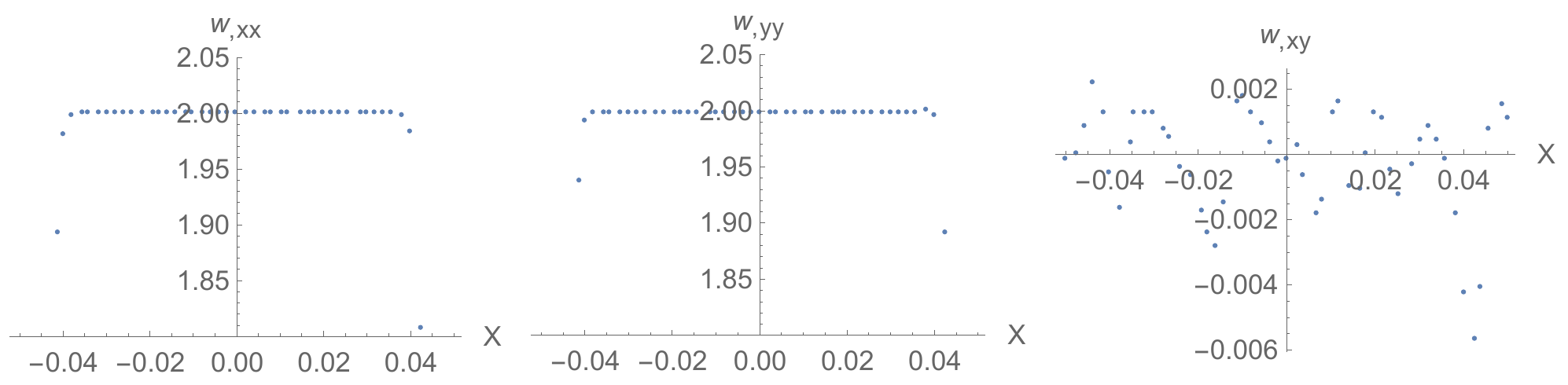}
\caption{The nonlocal Hessian with 240 neighbors in each support. }
\label{fig:N240HessianDeriInhomo}
\end{figure}

For regular particle distribution in 2D, small number of particles in support can accurately define the nonlocal Hessian operator, as shown in Fig.\ref{fig:N10HessianDeri}, Fig.\ref{fig:N24HessianDer} and Fig.\ref{fig:N48HessianDer}.
For irregular particle distribution as shown in Fig.\ref{fig:voronoiParticle}, larger number of particles in support are required to define the nonlocal Hessian operator, as shown in Fig.\ref{fig:N120HessianDeriInhomo} and Fig.\ref{fig:N240HessianDeriInhomo}.

\begin{figure}[htp]
\centering
\includegraphics[width=14cm]{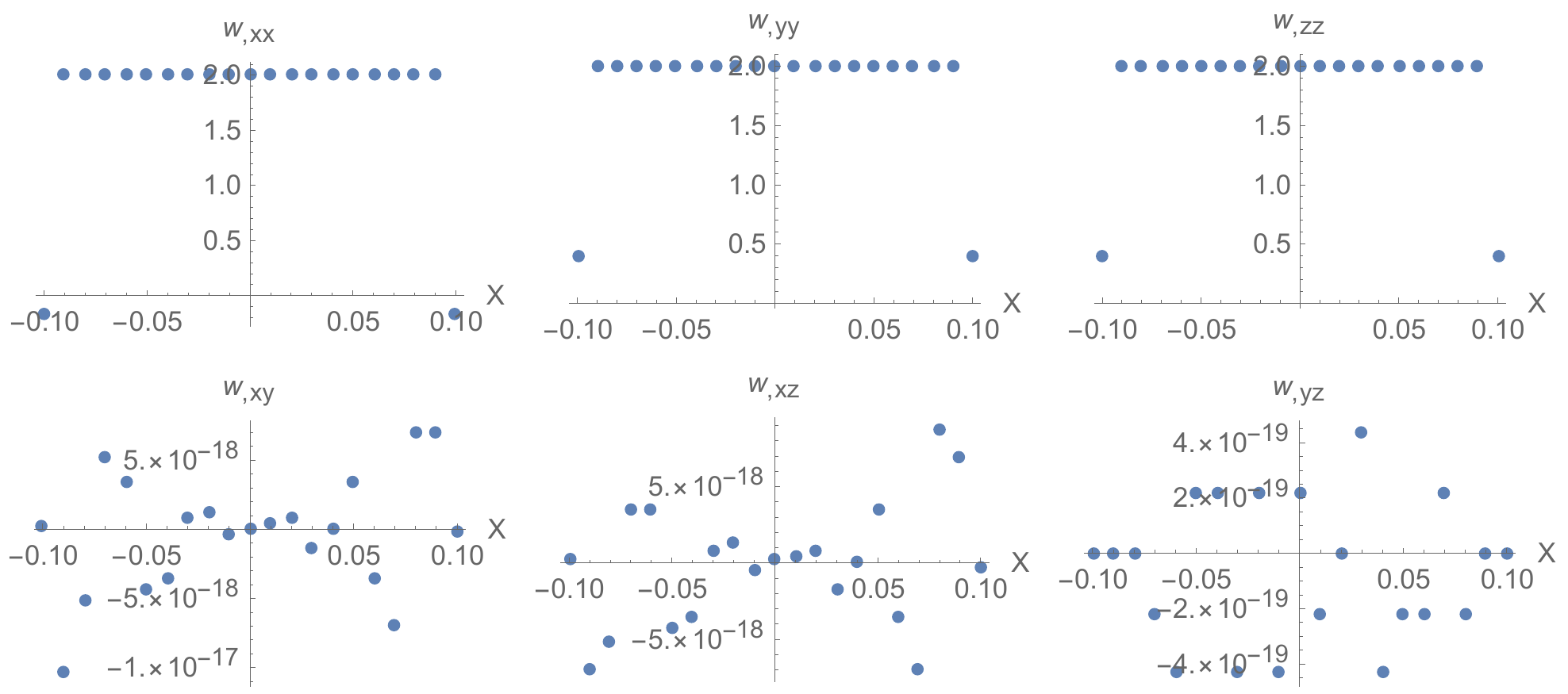}
\caption{The nonlocal Hessian with 26 neighbors in each support. }
\label{fig:N26HessianDer3D}
\end{figure}

\begin{figure}[htp]
\centering
\includegraphics[width=14cm]{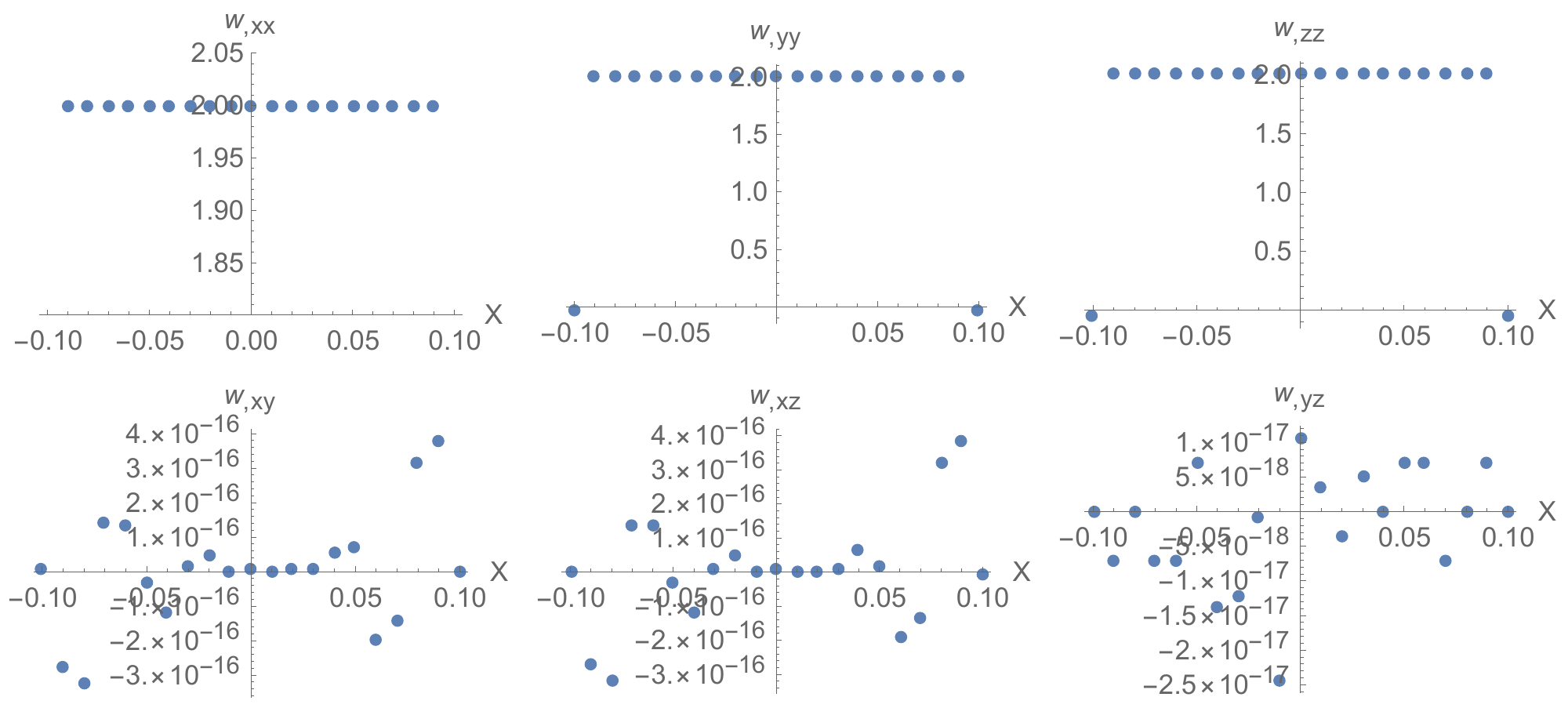}
\caption{The nonlocal Hessian with 48 neighbors in each support. }
\label{fig:N48HessianDer3D}
\end{figure}

\begin{figure}[htp]
\centering
\includegraphics[width=12cm]{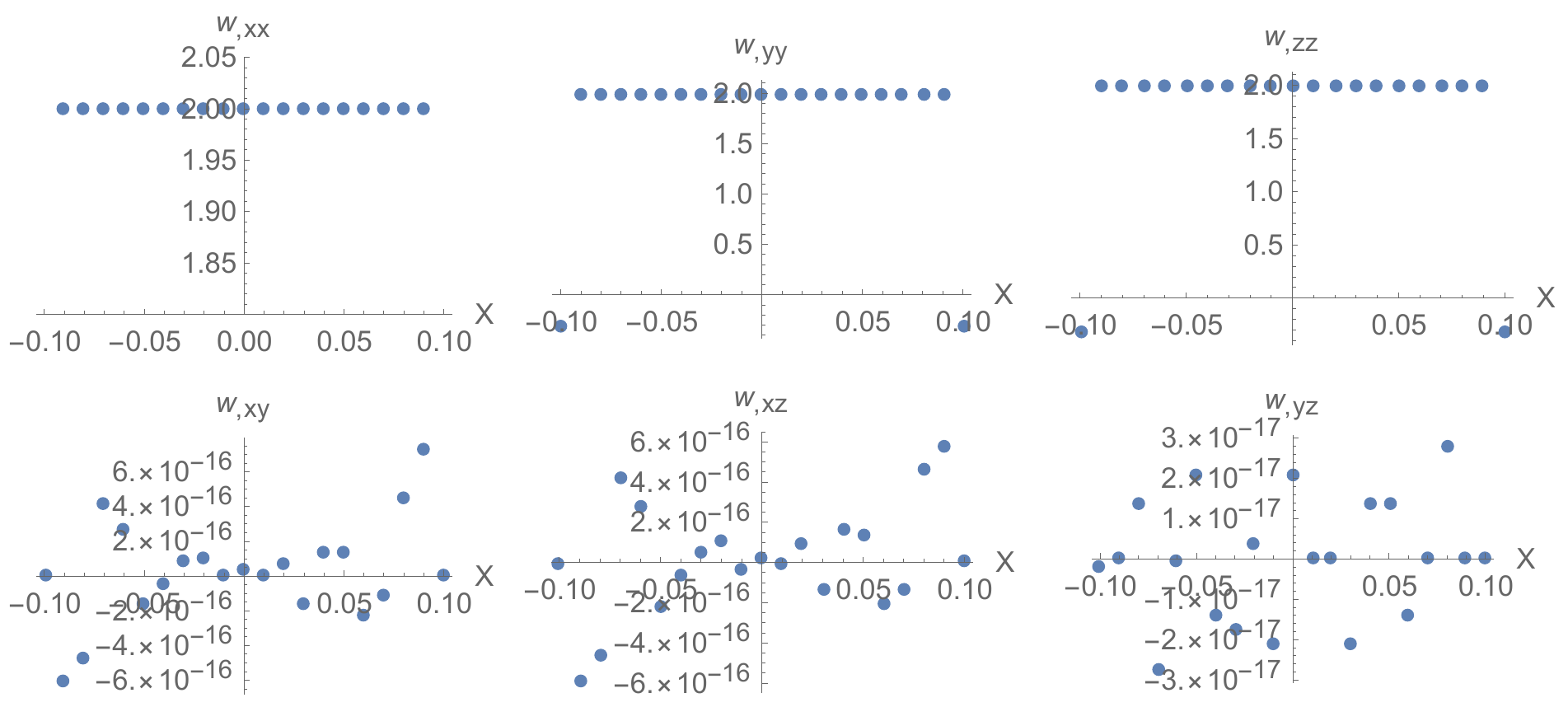}
\caption{The nonlocal Hessian with 64 neighbors in each support. }
\label{fig:N64HessianDer3D}
\end{figure}

For regular particle distribution in 3D, small number of particles in support can accurately define the nonlocal Hessian operator, as shown in Fig.\ref{fig:N26HessianDer3D}, Fig.\ref{fig:N48HessianDer3D} and Fig.\ref{fig:N64HessianDer3D}.

\subsection{Square thin plate subject to pressure}
\begin{figure}[htp]
\centering
\includegraphics[width=4cm]{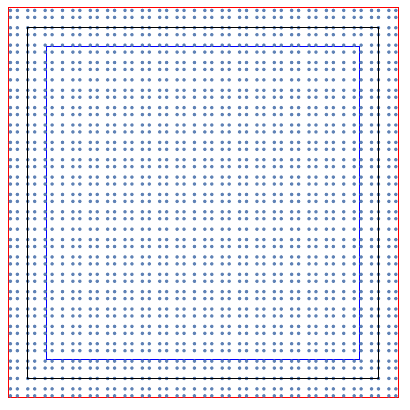}
\caption{The implementation of clamped boundary condition. The particles in black rectangle represent the physical model and particles outside of the blue rectangle are applied with penalty $p^{hg}=400 E$.}
\label{fig:ClampedPlateParticle}
\end{figure}

\begin{figure}[htp]
\centering
\subfigure[]{
\label{fig:s4td97}
\includegraphics[width=.4\textwidth]{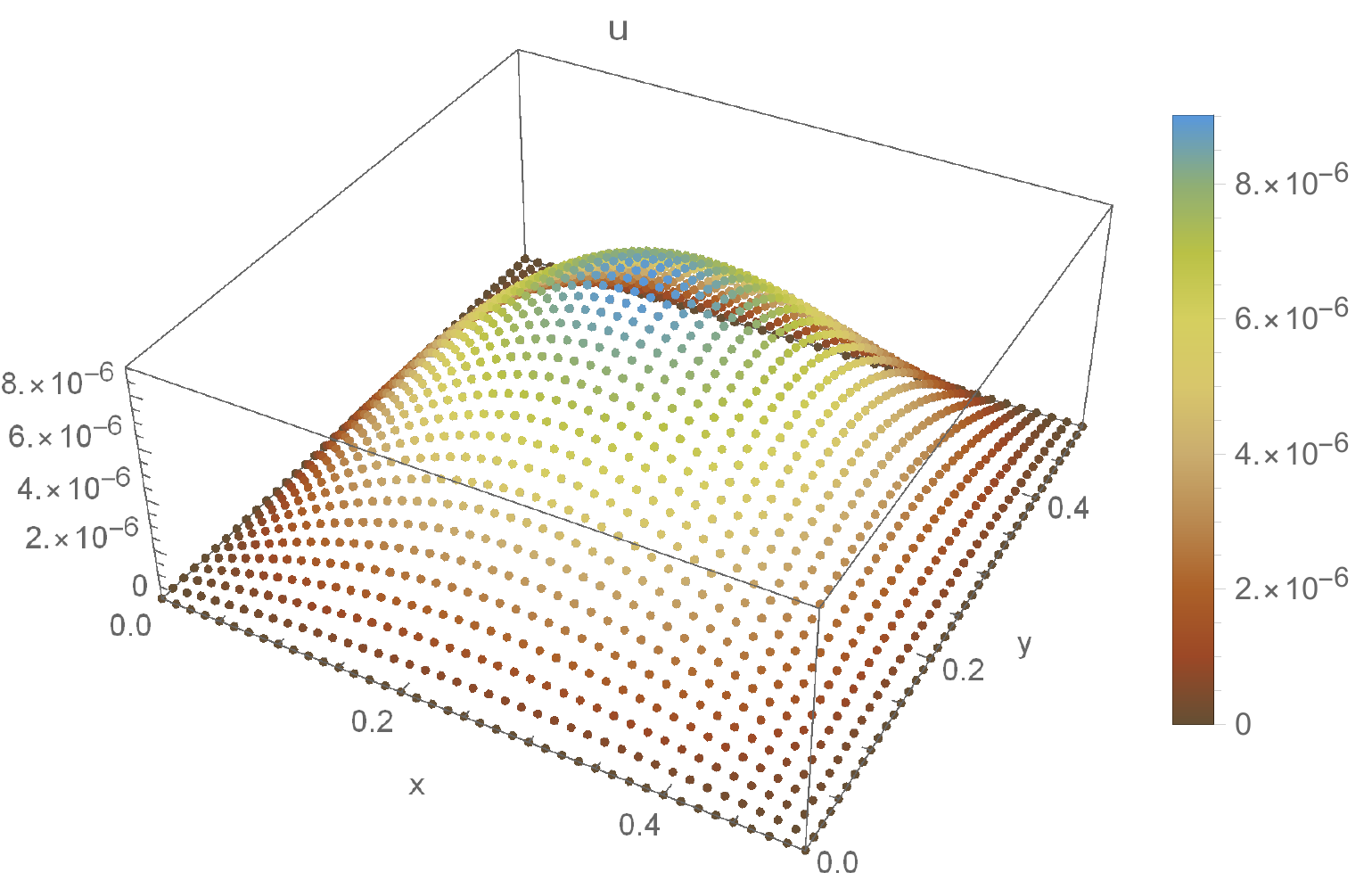}}
\vspace{.1in}
\subfigure[]{
\label{fig:s4t2d9}
\includegraphics[width=.4\textwidth]{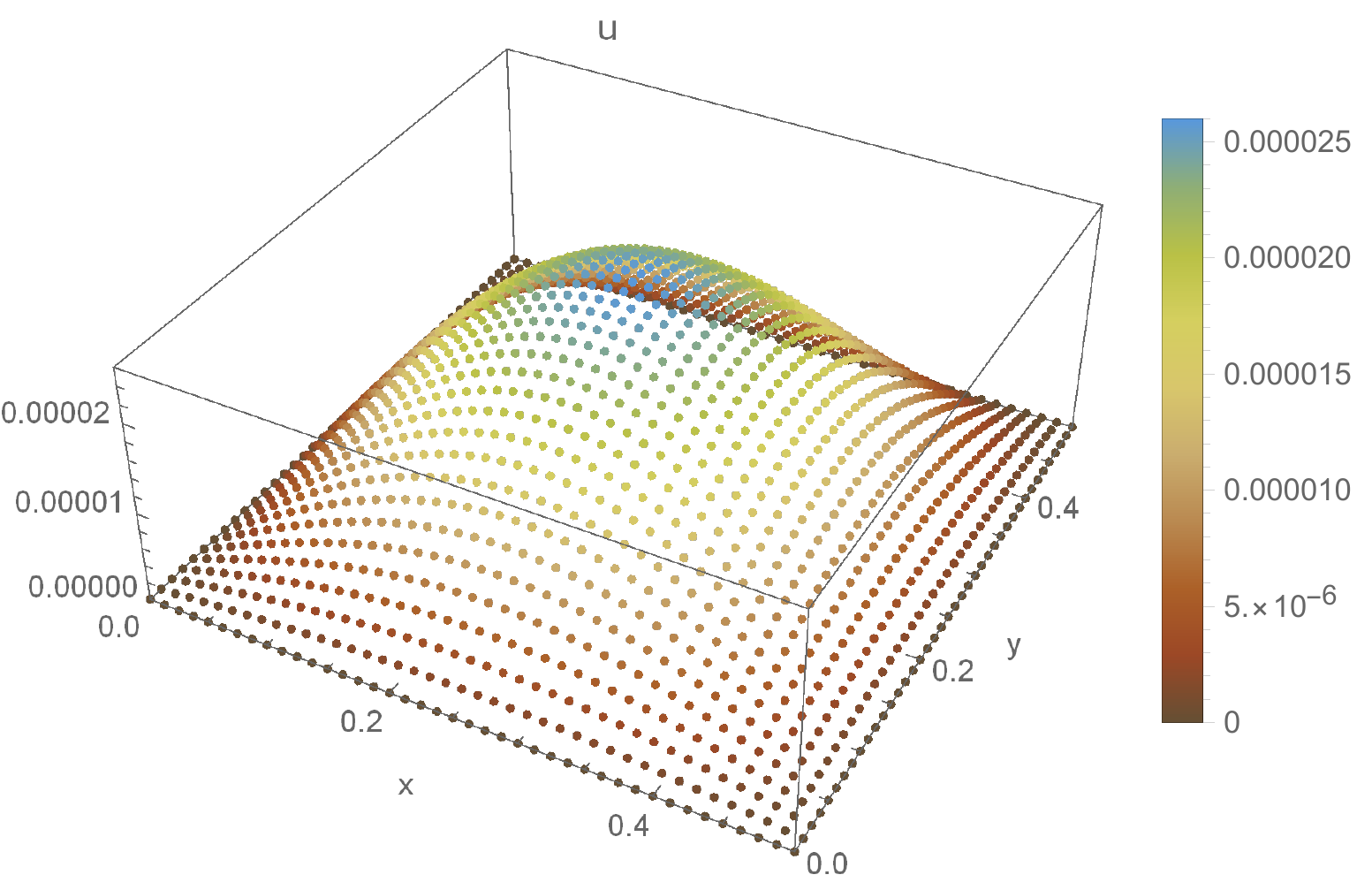}}\\
\subfigure[]{
\label{fig:s4t4d87}
\includegraphics[width=.4\textwidth]{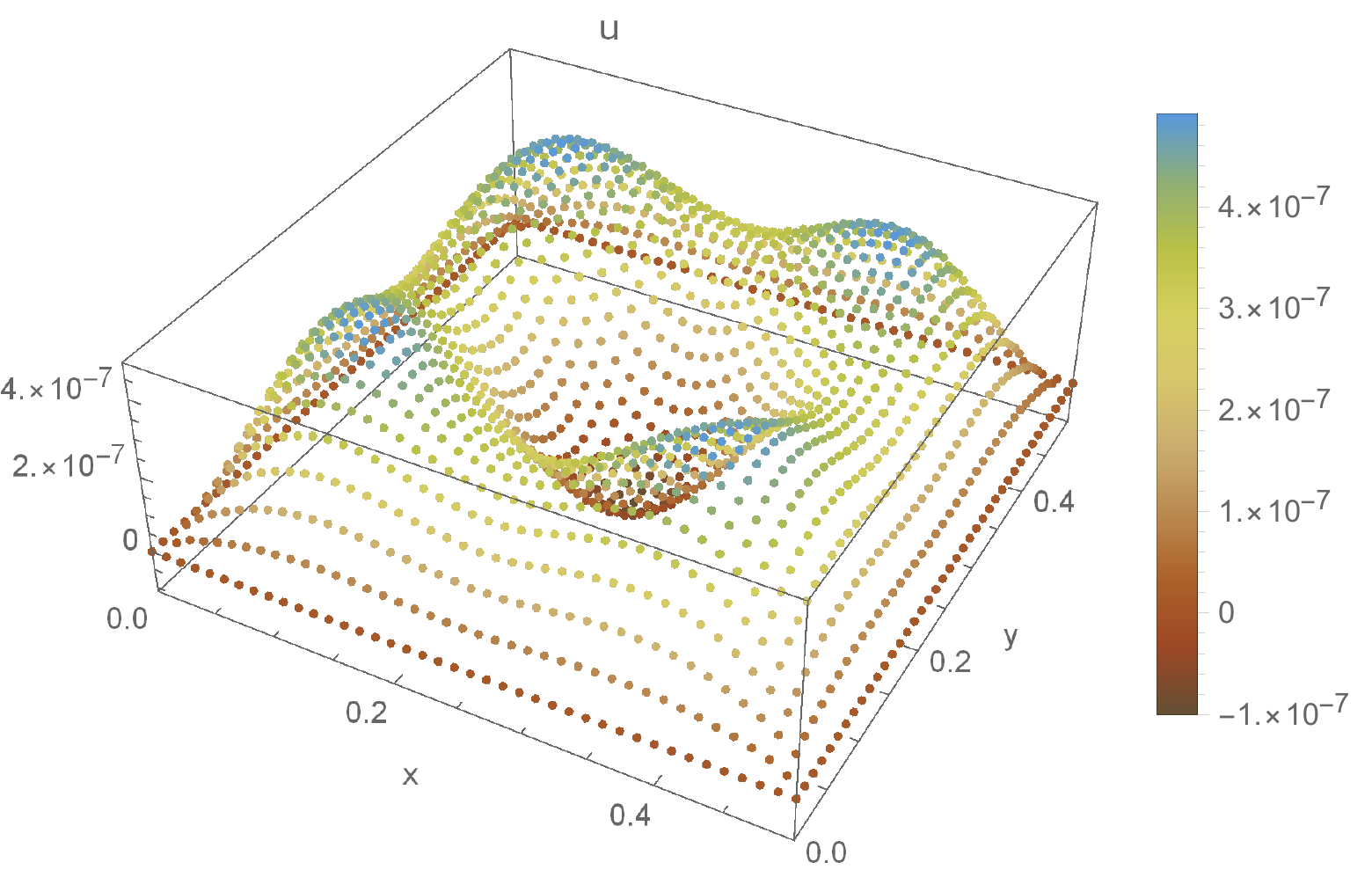}}
\vspace{.1in}
\subfigure[]{
\label{fig:s4t6d77}
\includegraphics[width=.4\textwidth]{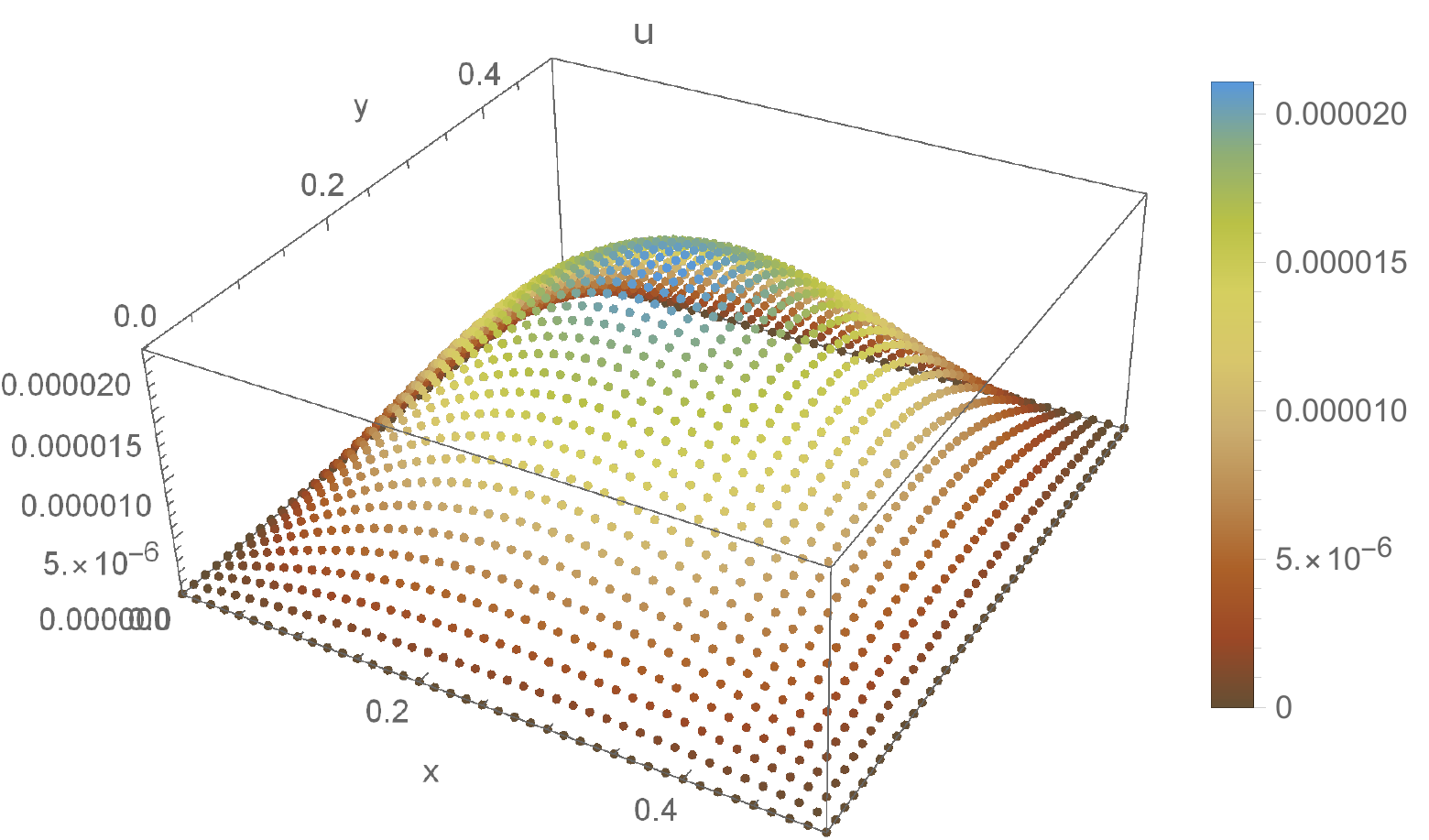}}\\
\vspace{.3in}
\caption{Deflection of simply supported plate at (a) $t=0.97$ ms (b) $t=2.9$ ms (c) $t=4.87$ ms and (d) $t=6.77$ ms.}\label{fig:simpleDeflection}
\end{figure}

\begin{figure}[htp]
\centering
\subfigure[]{
\label{fig:Clamp40t1}
\includegraphics[width=.4\textwidth]{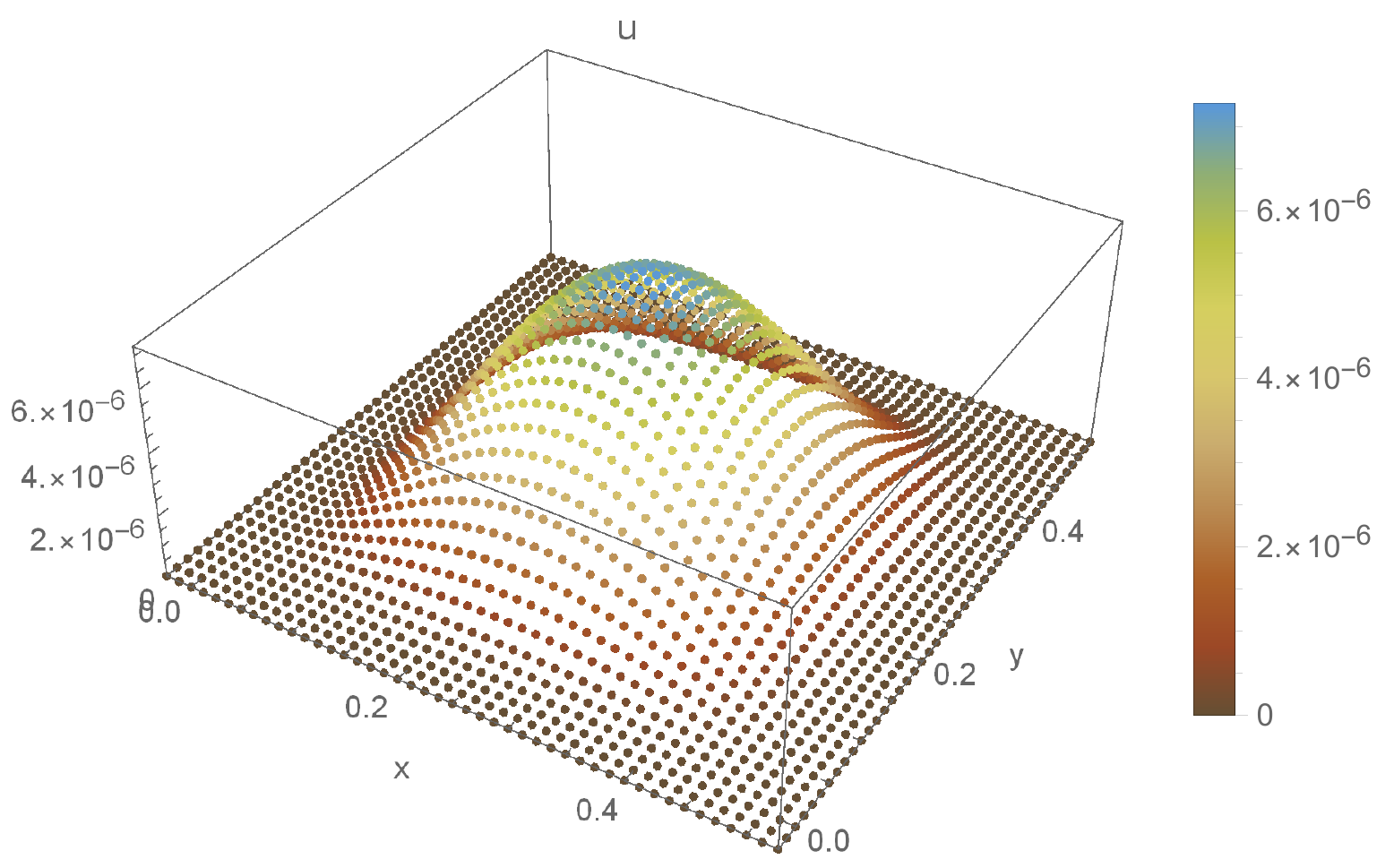}}
\vspace{.1in}
\subfigure[]{
\label{fig:Clamp40t2}
\includegraphics[width=.4\textwidth]{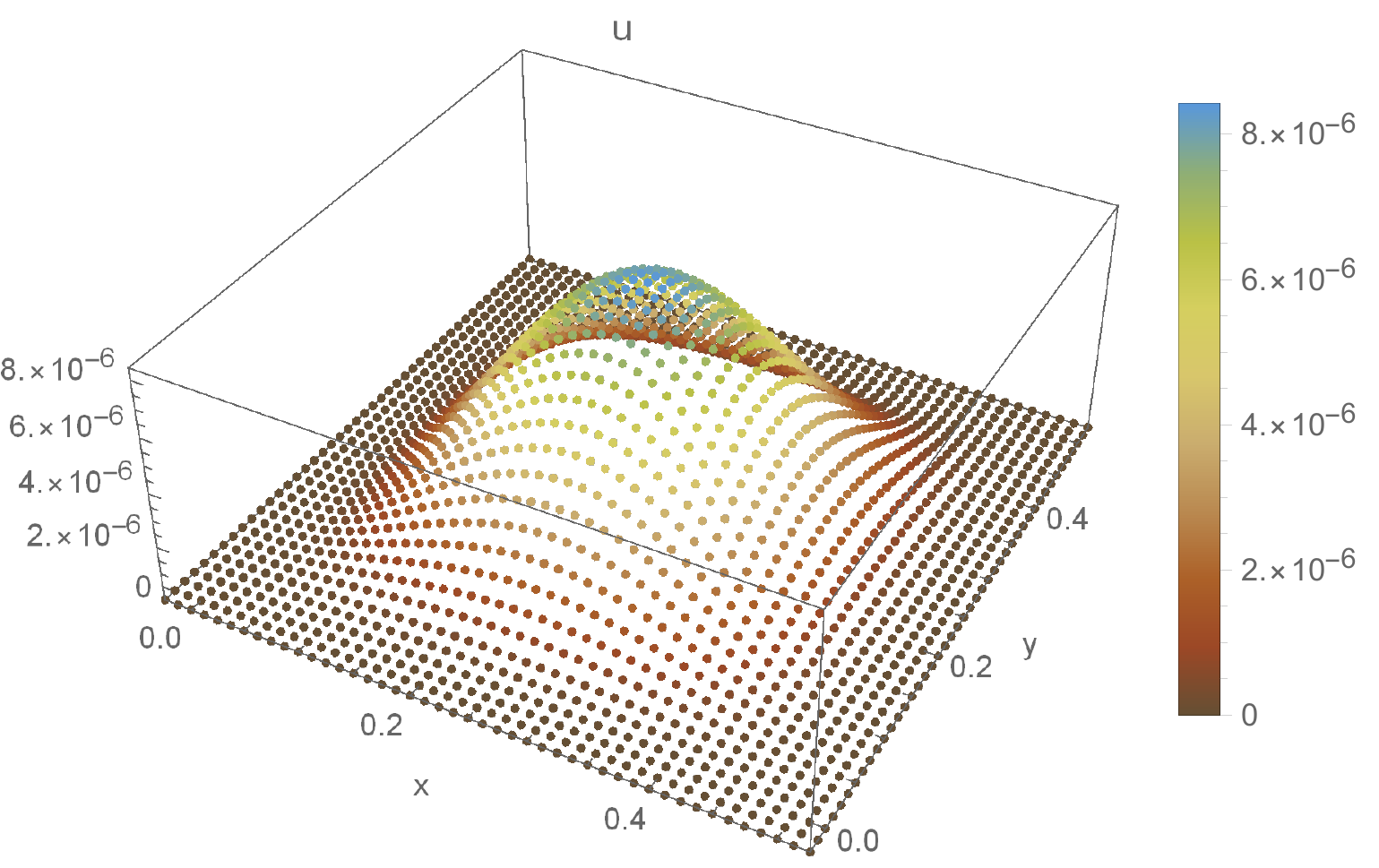}}\\
\subfigure[]{
\label{fig:Clamp40t3}
\includegraphics[width=.4\textwidth]{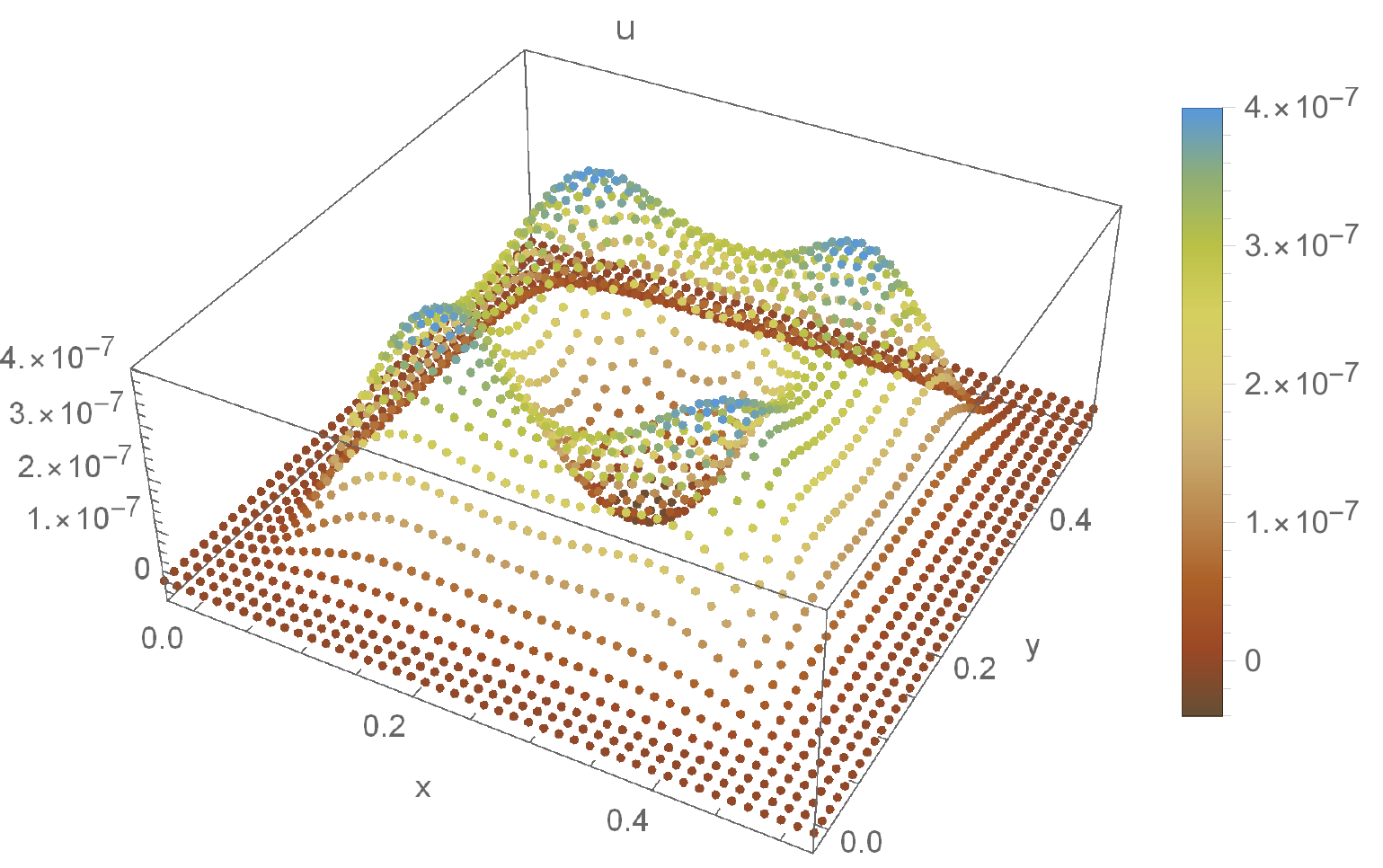}}
\vspace{.1in}
\subfigure[]{
\label{fig:Clamp40t4}
\includegraphics[width=.4\textwidth]{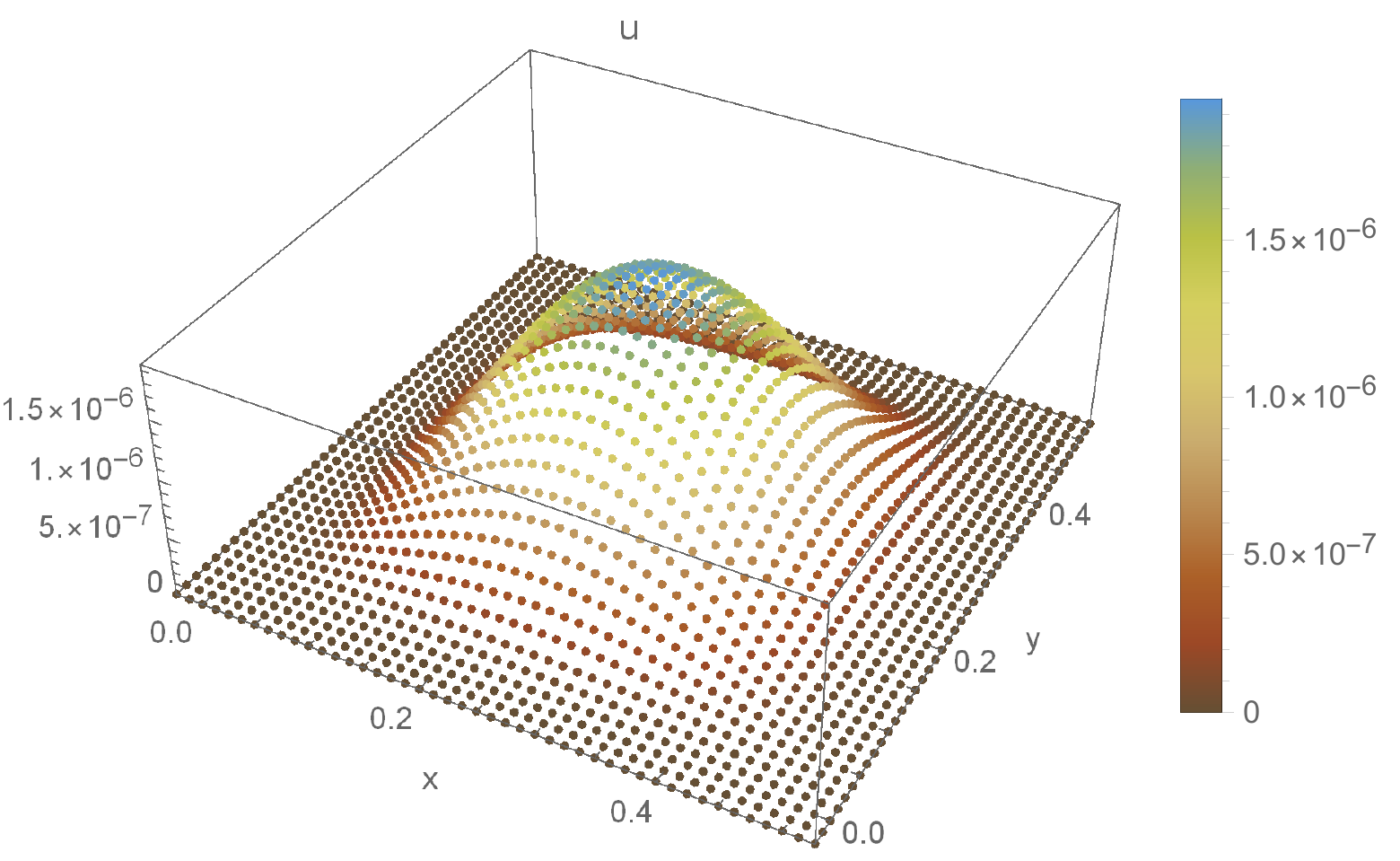}}\\
\vspace{.3in}
\caption{Deflection of simply supported plate at (a) $t=0.966$ ms (b) $t=1.44$ ms (c) $t=2.42$ ms and (d) $t=2.90$ ms.}\label{fig:ClampDeflection}
\end{figure}

\begin{figure}[htp]
\centering
\subfigure[]{
\label{fig:Simple40Deflection}
\includegraphics[width=.4\textwidth]{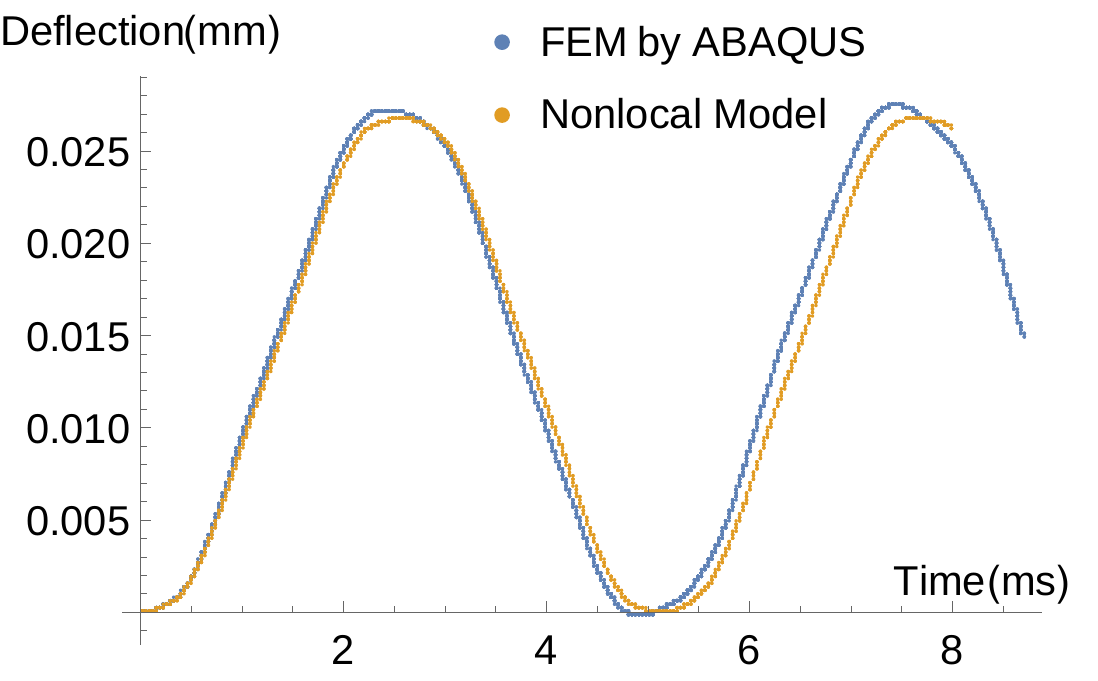}}
\vspace{.1in}
\subfigure[]{
\label{fig:Clamp40Deflection}
\includegraphics[width=.45\textwidth]{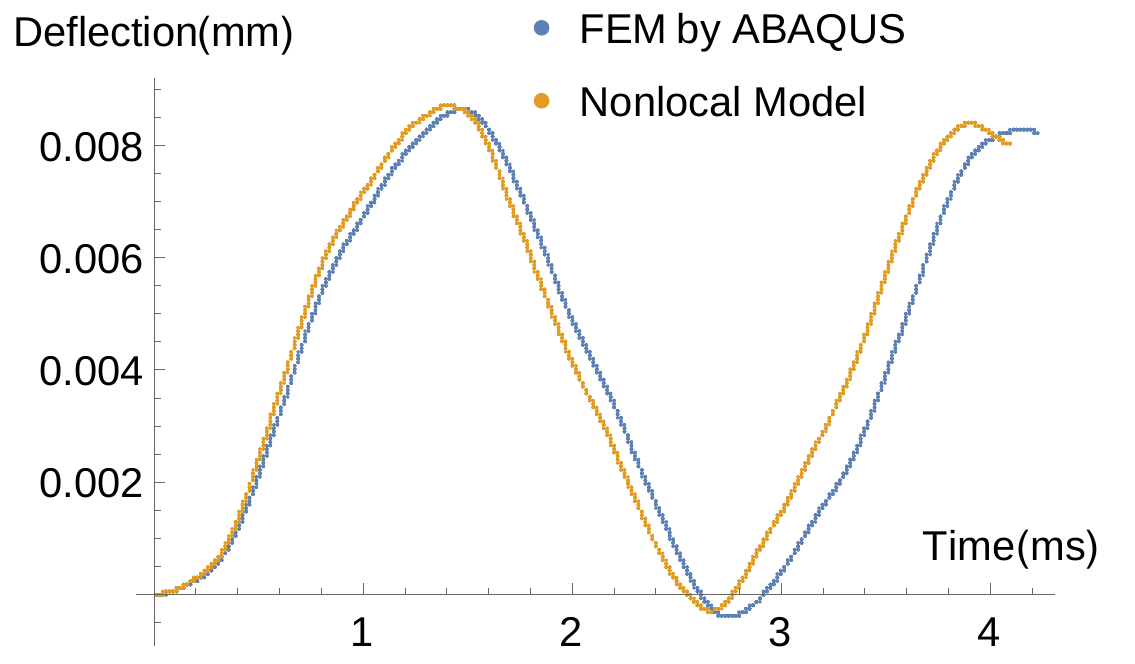}}\\
\caption{Deflection of central point for (a) simply support plate and (b) clamped support plate.}
\end{figure}

The dimensions of the plate are $0.5\times 0.5$ $m^2$ with a thickness of $0.01$ m. The material parameters are elastic modulus $E=210$ GPa, Poisson ratio $\nu=0.3$. The plate is applied with a static pressure load of $p=10^3$ Pa. Two boundary conditions are taken into account: a) four sides are all simply supported and b) four sides are all clamped. The case of clamped boundary constrains the rotation as well as the deflection. The reference result is calculated by $64\times 64$ S4R elements in ABAQUS without considering the geometric nonlinearity. For the simply supported boundary conditions, the particles on the boundaries of the plate are fixed. The enforcement of clamped boundary conditions requires some special treatment. As shown Fig.\ref{fig:ClampedPlateParticle}, the particles in the black rectangle are the particles in the physical model and the particles outside of the blue rectangle are applied with penalty $p^{hg}=400 E$ while the particles inside the blue rectangle with penalty $p^{hg}=0$. The deflection for a simply supported plate at different times are plotted in Fig.\ref{fig:simpleDeflection}. The deflections for a clamped plate at different times are depicted in Fig.\ref{fig:ClampDeflection}. The deflection of the central point of the plate is monitored and compared with the result by ABAQUS, as shown in Fig.\ref{fig:Simple40Deflection} and Fig.\ref{fig:Clamp40Deflection}, where good agreement with FEM model is observed.

%s4deflectiond97.pdf,s4deflection2d9.pdf,s4deflection4d87.pdf,s4deflection6d77.pdf

%{0.96574, 1.44476, 2.41617, 2.89917}

%\begin{figure}[htp]
%\centering
%\subfigure[]{
%\label{fig:Deformation3gS10}
%\includegraphics[width=.45\textwidth]{Deformation3gS10.pdf}}
%\vspace{.1in}
%\subfigure[]{
%\label{fig:ContourOfEnergyAll3}
%\includegraphics[width=.4\textwidth]{ContourOfEnergyAll3.pdf}}\\
%\vspace{.3in}
%\caption{(a) Deformation in $x$-direction, scaled by 10 times and (b) the distribution of total strain energy density for discretization of $120^2$ particles.}
%\end{figure}
%\subsection{Fracture modelling in elasticity}

\subsection{Single-edge notched tension test}
\begin{figure}[htp]
\centering
\includegraphics[width=4cm]{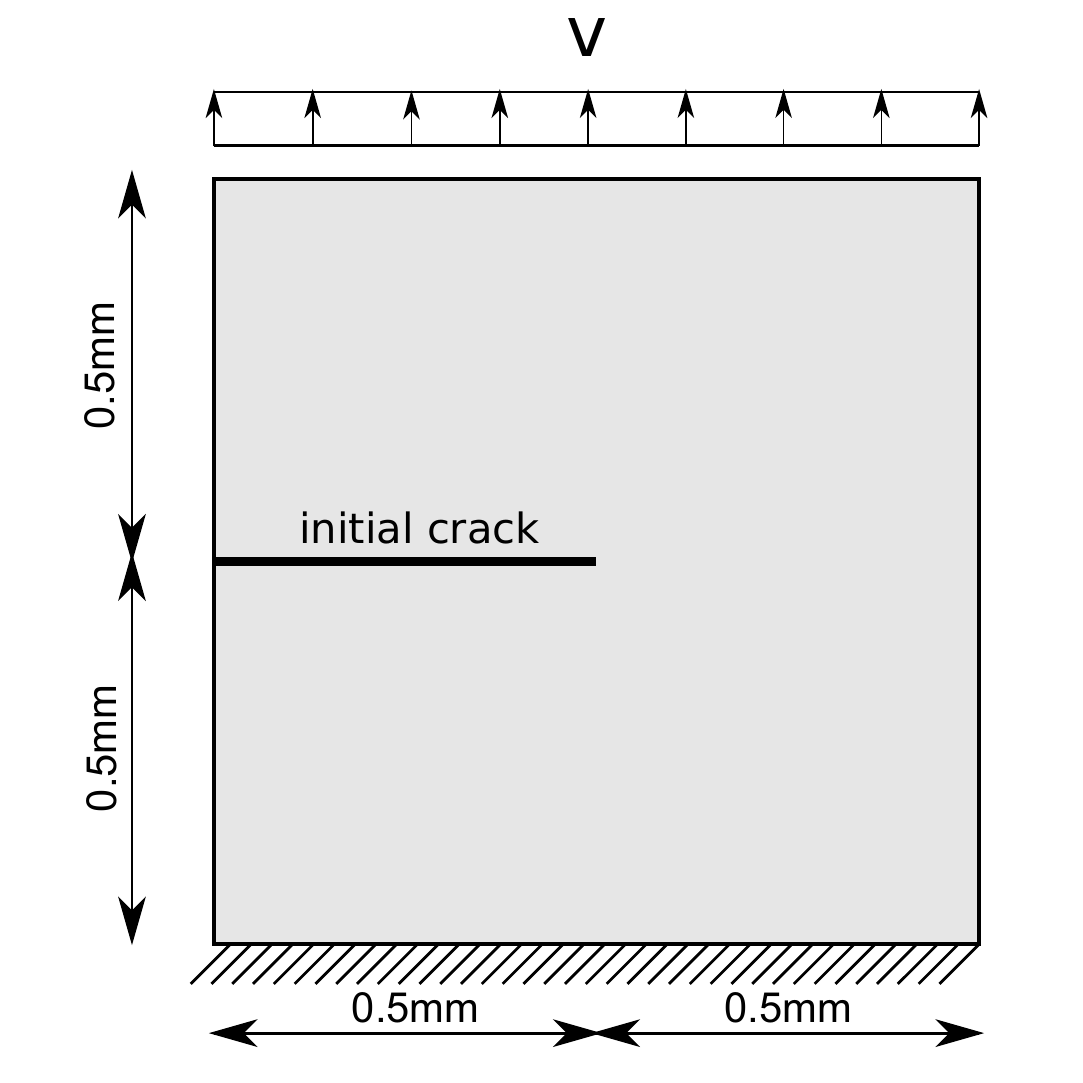}
\caption{Setup of the plate.}
\label{fig:platetensionGeom}
\end{figure}

\begin{figure}[htp]
\centering
\subfigure[]{
\label{fig:uy2d16}
\includegraphics[width=.45\textwidth]{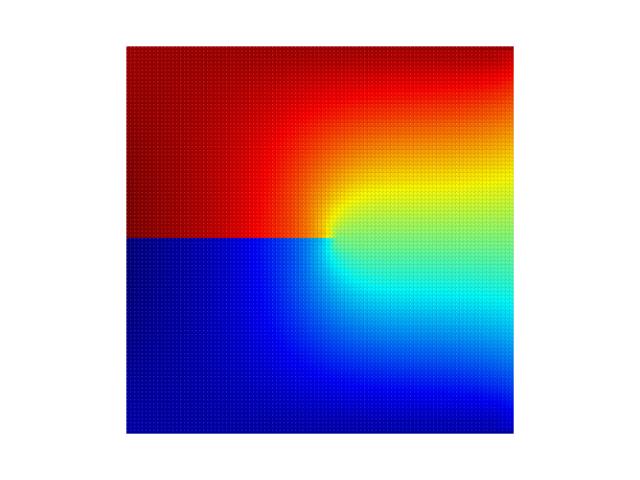}}
\vspace{.1in}
\subfigure[]{
\label{fig:uy2d19}
\includegraphics[width=.45\textwidth]{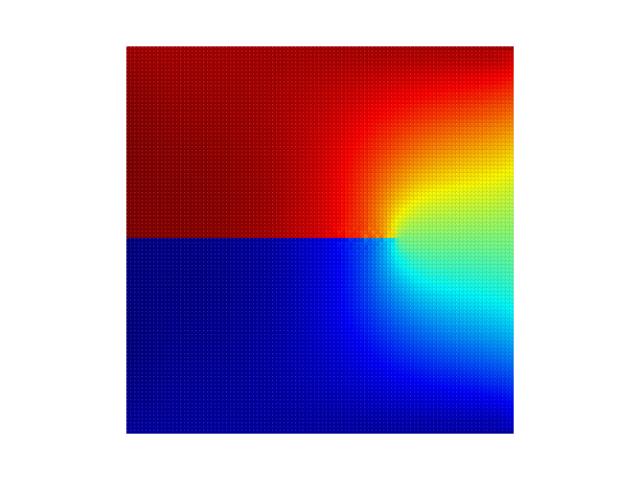}}\\
\caption{Displacement field at (a) $u_y=5.0 \times 10^{-3}$ mm, and (b) $u_y=5.5 \times 10^{-3}$ mm.}
\end{figure}

In this example, we tested the nonlocal elasticity by Eq.\ref{eq:dsnome} for single-edge notched tension in 2D under plane stress condition. The geometry setup is given in Fig.\ref{fig:platetensionGeom}. The bottom is fixed while the top of the plate is applied with velocity boundary condition $v=1$ m/s, which can achieve the quasi-static condition. The material parameters are $E=210\mbox{ GPa, }\nu=0.3$ and critical strain is set as $s_{max}=0.02$. The plate is discretized into 100$\times$100 particles. Each particle's support consists of 33 nearest neighbours. The initial crack is created by modifying the neighbour list when searching the nearest neighbours. The fixed number of neighbours in support results in particles near the boundary with relatively large support sizes and particles in the centre of the plate with small support sizes. A duration of $T=6.5\times 10^{-6}$ seconds is integrated by approximately 4500 steps at a time increment of $\Delta t=1.5418\times 10^{-9}$ seconds. The displacement field $u_y$ at step 3250 and step 4200 are depicted in Fig.\ref{fig:uy2d16} and Fig.\ref{fig:uy2d19}, respectively.

\begin{figure}[htp]
\centering
\subfigure[]{
\label{fig:UyField2}
\includegraphics[width=.41\textwidth]{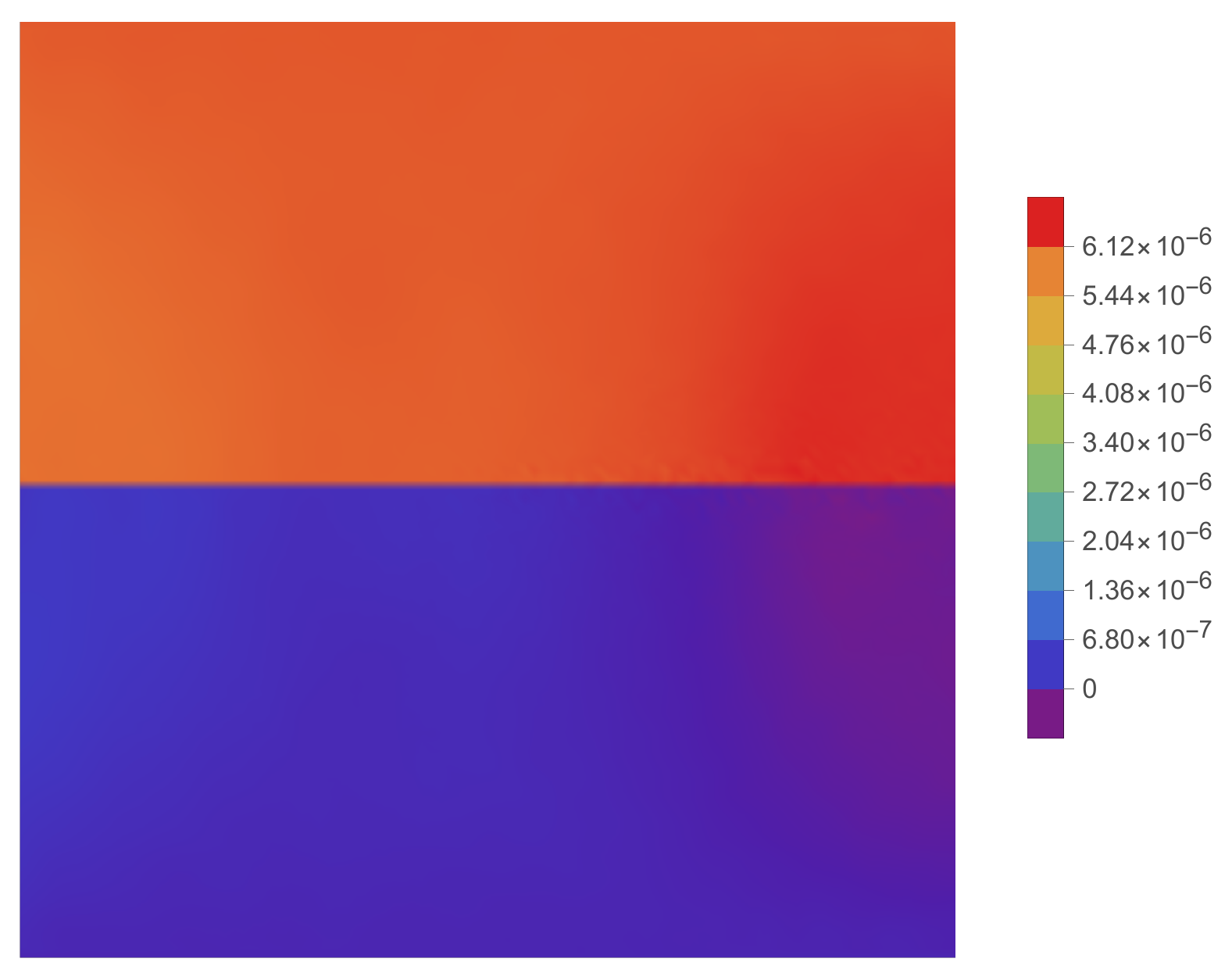}}
\vspace{.1in}
\subfigure[]{
\label{fig:hourglassEnergy1}
\includegraphics[width=.4\textwidth]{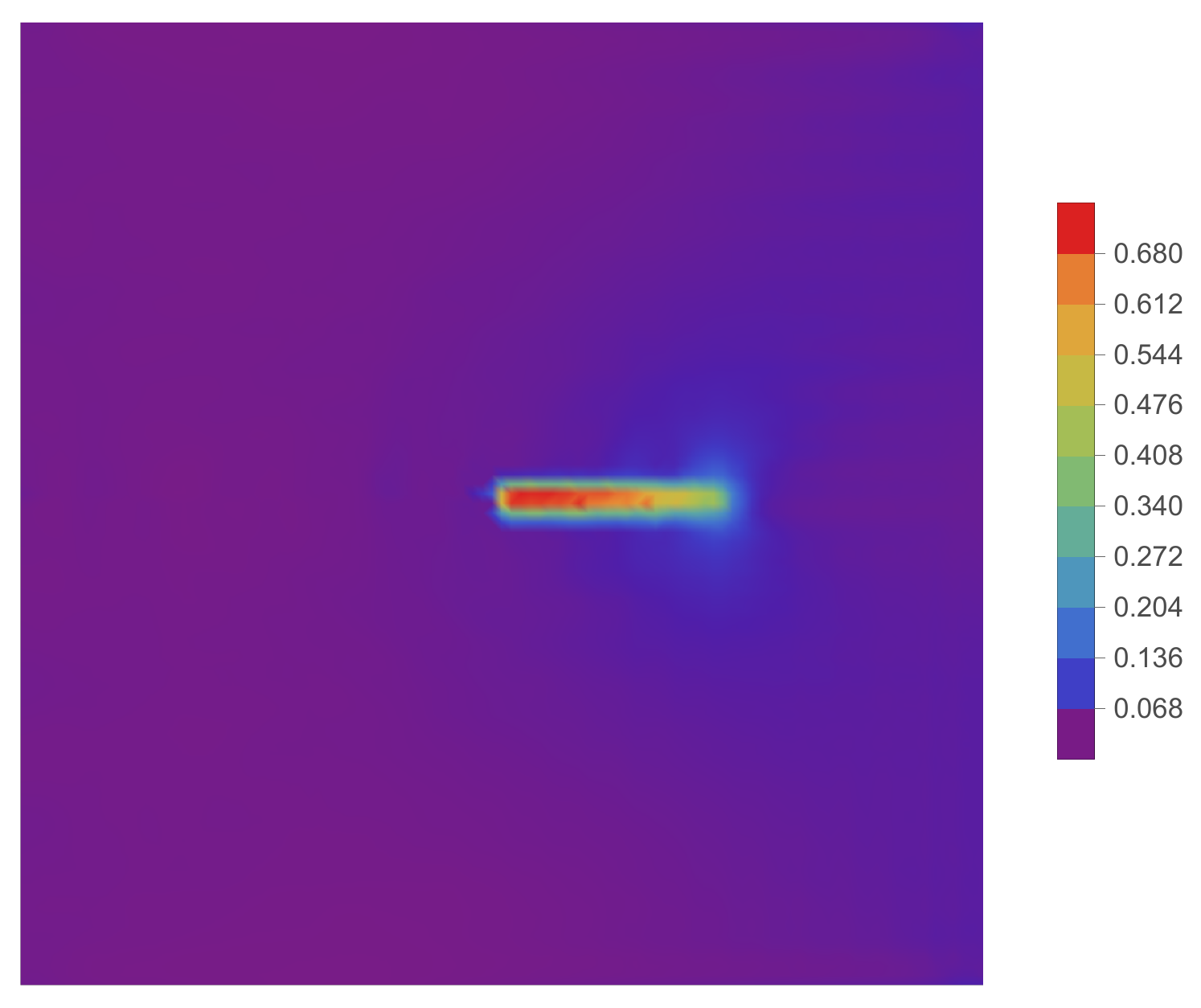}}\\
\subfigure[]{
\label{fig:damage1}
\includegraphics[width=.4\textwidth]{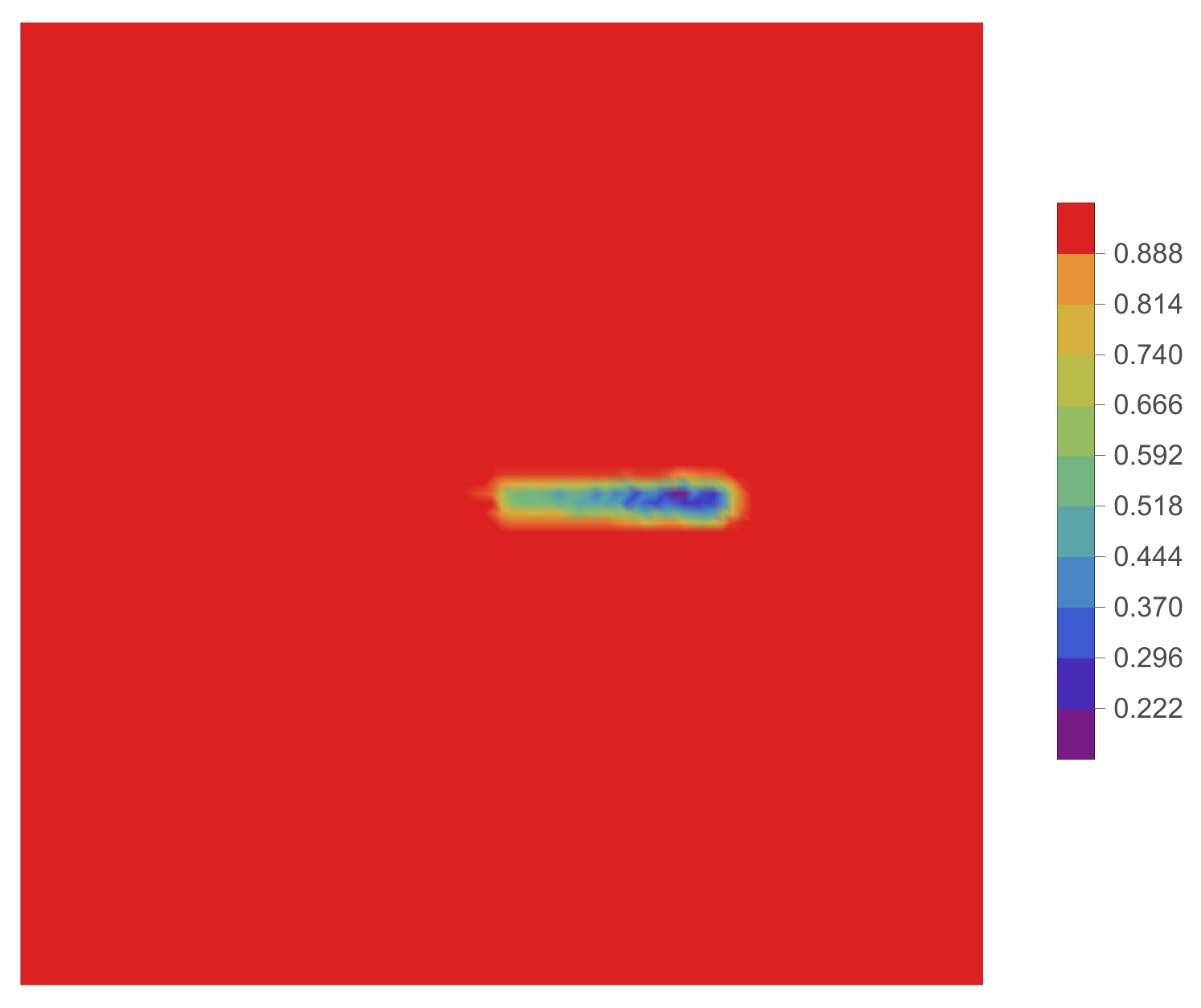}}
\vspace{.1in}
\subfigure[]{
\label{fig:damage2}
\includegraphics[width=.4\textwidth]{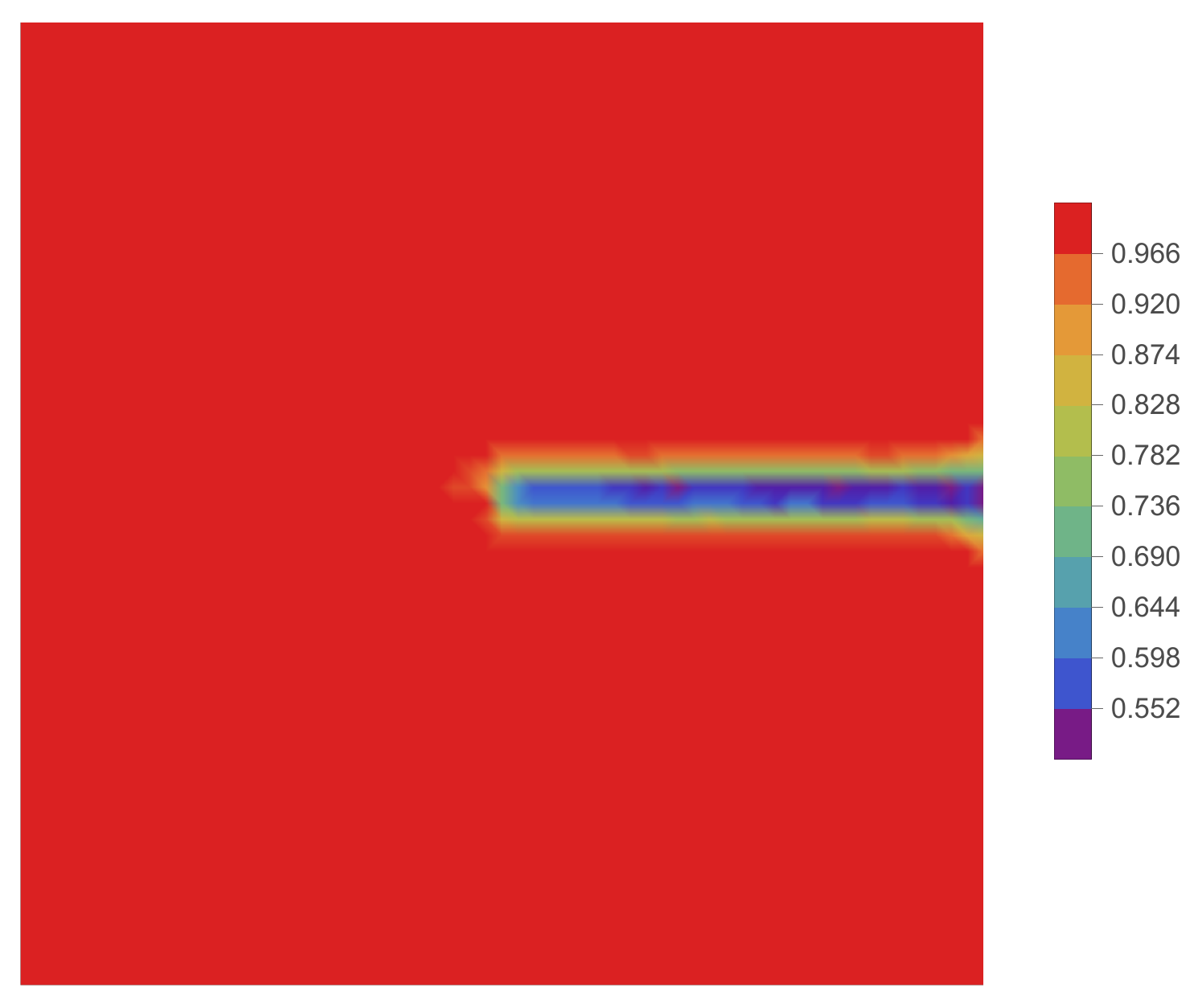}}\\
\vspace{.3in}
\caption{(a) Displacement $u_y$ at full damage (b) damage field at $u_y=5.5 \times 10^{-3}$ mm (c) damage field at $u_y=6.2 \times 10^{-3}$ mm and (d) operator energy at $u_y=5.5 \times 10^{-3}$ mm.}\label{fig:plateDamage}
\end{figure}%{0.96574, 1.44476, 2.41617, 2.89917}
Fig.\ref{fig:UyField2} is the displacement field $u_y$ at full damage, where the interaction of internal force between the two half planes is cut and rigid body displacement dominates. Fig.\ref{fig:hourglassEnergy1} is the distribution of hourglass energy. We can observe that the hourglass energy is concentrated on the crack surface and crack front tip. Fig.\ref{fig:damage1} and Fig.\ref{fig:damage2} are the snapshots of damage field, which confirms that the instability criterion in \S \ref{sec:fracCriterion} is stable for fracture modelling.

\begin{figure}[htp]
\centering
\subfigure[]{
\label{fig:energyCurve}
\includegraphics[width=.45\textwidth]{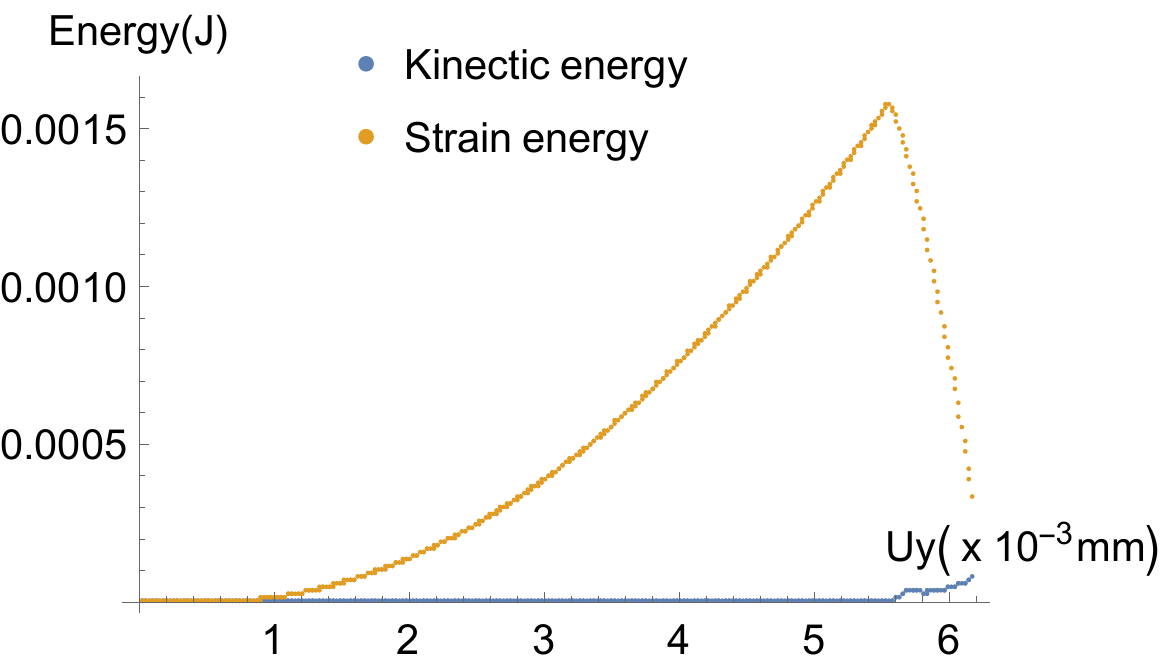}}
\subfigure[]{
\label{fig:loadCurve}
\includegraphics[width=.45\textwidth]{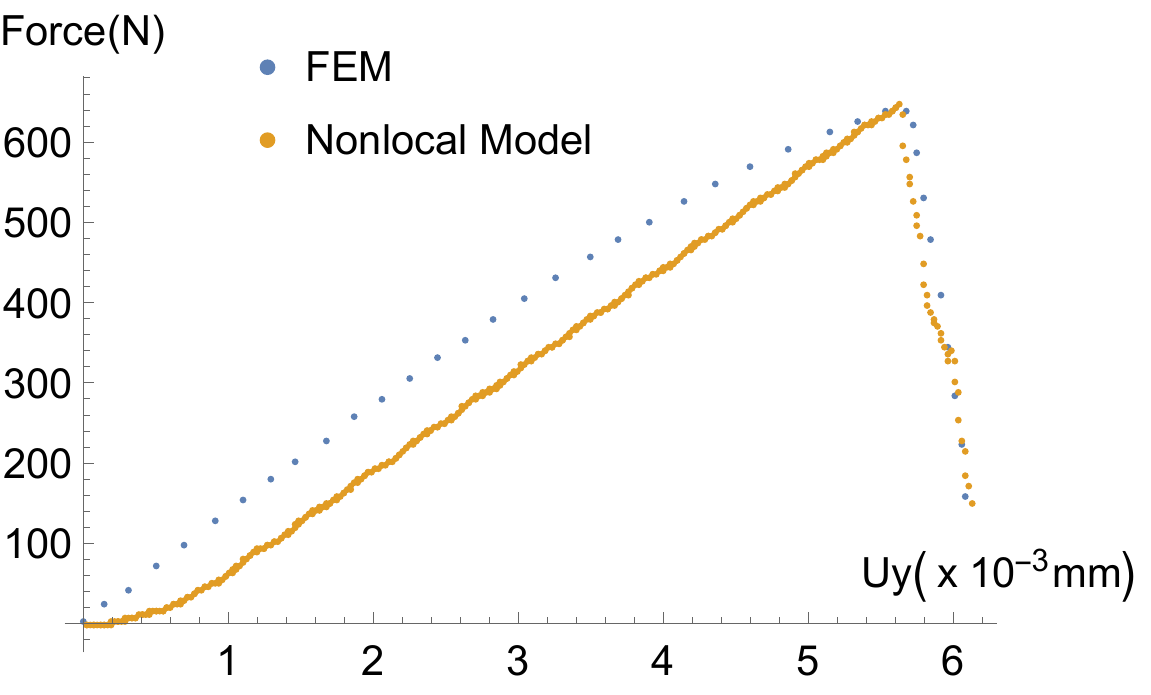}}\\
\caption{(a) Load curve on displacement; (b) energy curve on displacement.}
\end{figure}

Although the plate is solved by an explicit dynamic method, the kinetic energy is much lower than the strain energy as shown by Fig.\ref{fig:energyCurve}. The dynamic load curve agrees well with that by the finite element method in Ref \cite{miehe2010thermodynamically}, as shown by Fig.\ref{fig:loadCurve}. One possible reason for the difference in reaction force increment is due to explicit algorithm and nonlocal effect of current formulation.

\subsection{Out-of-plane shear fracture in 3D}
\begin{figure}[htp]
\centering
\includegraphics[width=5cm]{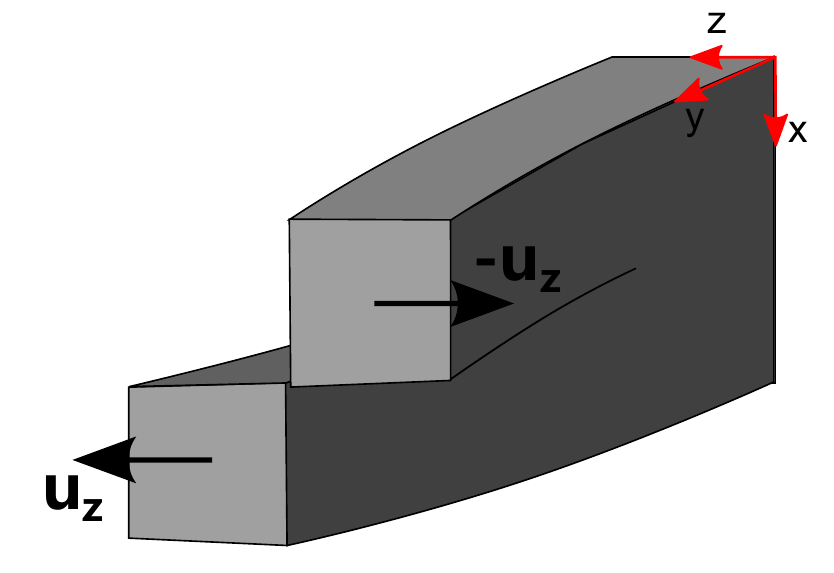}
\caption{Illustration of out-of-plane shear fracture.}
\label{fig:mode3f}
\end{figure}

\begin{figure}[htp]
\centering
\includegraphics[width=8cm]{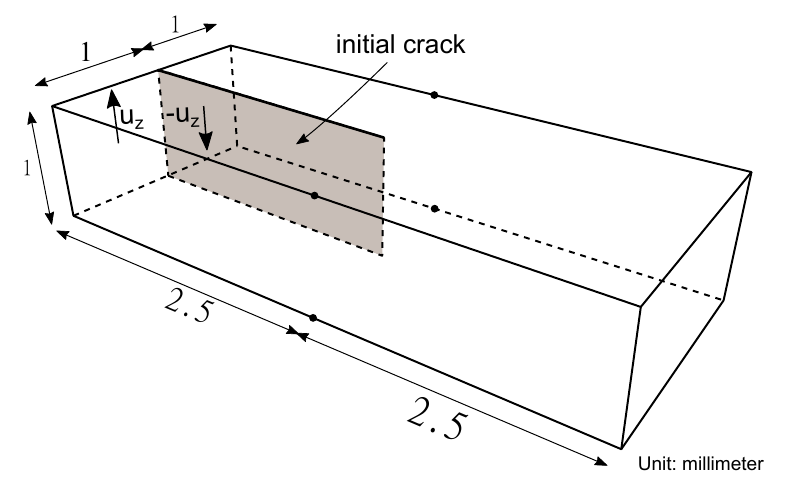}
\caption{Setup of the specimen.}
\label{fig:mode3f2}
\end{figure}

\begin{figure}[htp]
\centering
\subfigure[]{
\label{fig:crack26}
\includegraphics[width=.4\textwidth]{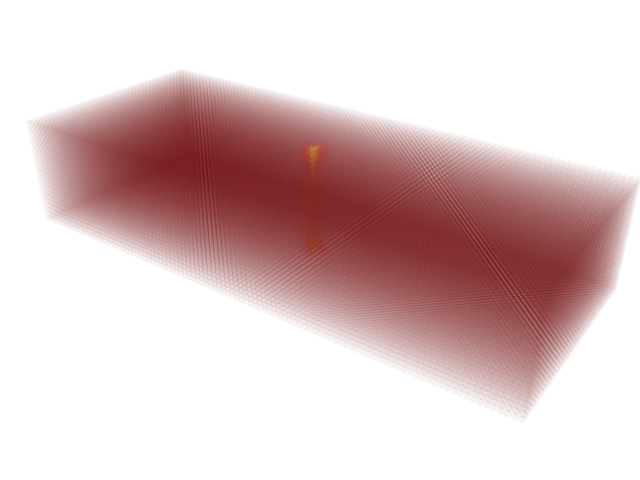}}
\vspace{.1in}
\subfigure[]{
\label{fig:crack35}
\includegraphics[width=.4\textwidth]{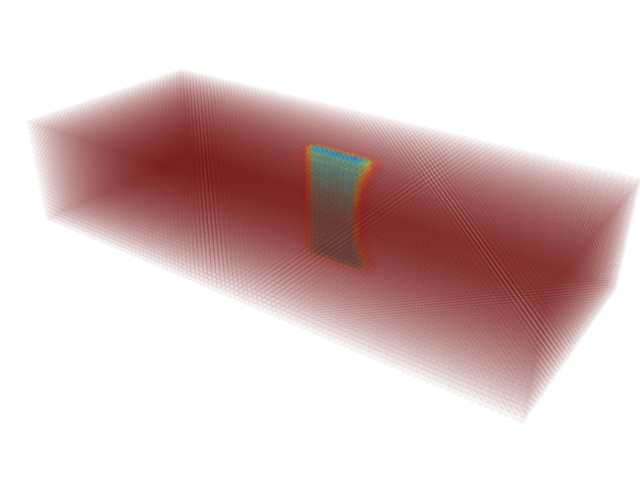}}\\
\subfigure[]{
\label{fig:crack50}
\includegraphics[width=.4\textwidth]{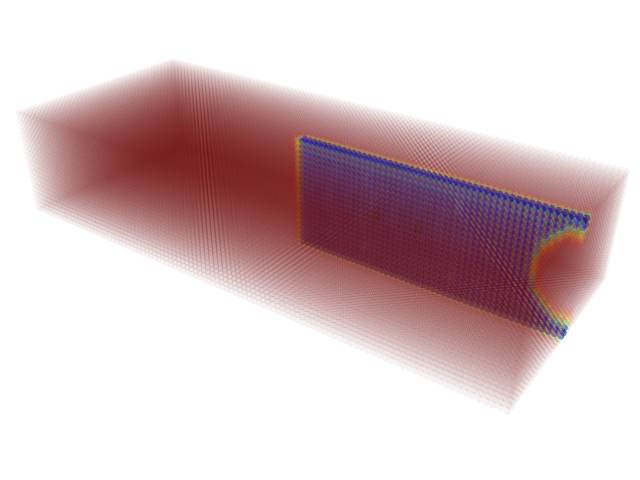}}
\vspace{.1in}
\subfigure[]{
\label{fig:crack51}
\includegraphics[width=.4\textwidth]{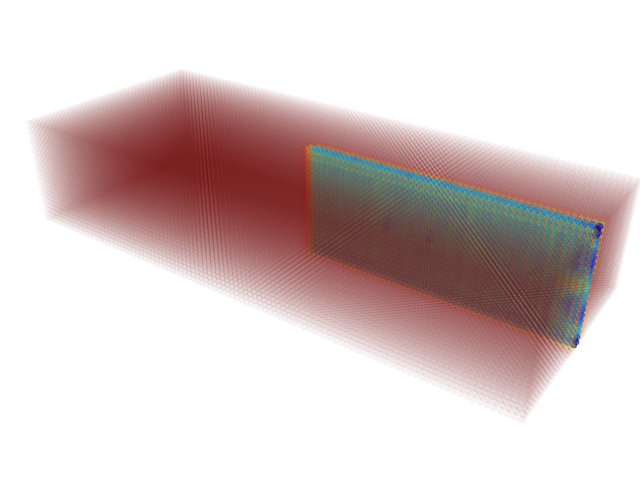}}\\
\vspace{.3in}
\caption{Crack surfaces at (a) step 1550 (b) step 2050 (c) step 2950 and (d) step 3000.}\label{fig:plateDamage3d}
\end{figure}

For brittle fracture, the basic modes of fracture are tensile fracture, in-plane shear fracture, out-of-plane shear fracture. In this section, we apply the instability damage criterion to the out-of-plane shear fracture, as shown in Fig.\ref{fig:mode3f}. The dimensions of the specimen are $5\times 2\times 1$ mm$^{3}$, as shown in Fig.\ref{fig:mode3f2}. The size of the initial crack surface is $2\times 1$ mm$^{2}$. The velocity boundary conditions $u_z=1$ m/s are applied. The model is discretized into 86961 particles with particle size $\Delta x=0.05$ mm. Each particle has 102 neighbours in its support. Material parameters include elastic modulus $E = 210\times 10^9$ Pa and Poisson ratio $\nu = 0.3$ and density $\rho = 7800$ kg/m$^3$. The time step is selected as $\Delta t=7.7\times 10^{-9}$ seconds. A total of 3000 steps are calculated. The crack surface starts to propagate at step 1550. The crack surface at different steps are depicted in Fig.\ref{fig:plateDamage3d}.

%{0.96574, 1.44476, 2.41617, 2.89917}

%
%
%\begin{figure}[htp]
% \centering
% \includegraphics[width=8cm]{uz51.png}
% \caption{Displacement field $u_z$.}
% \label{fig:platetensionGeom}
%\end{figure}

\section{Conclusion}
In this paper, we employ the recent proposed NOM to derive the nonlocal strong forms for various physical models, including elasticity, thin plate, gradient elasticity, electro-magneto-elastic coupled model and phase field fracture model. These models require second order partial derivative at most and we make use of the second-order NOM scheme, which contains the nonlocal gradient and nonlocal Hessian operator. Considering the fact that most physical models are compatible with the variational principle/weighted residual method, we start from the energy form/weak form of the problem, by inserting the nonlocal expression of the gradient/Hessian operator into the weak form, based on the dual property of the dual-support in NOM, the nonlocal strong form is obtained with ease. Such a process can be extended to many other physical problems in other fields. The derived strong forms are variationally consistent and allow elegant description for inhomogeneous nonlocality in both theoretical derivation and numerical implementation.

We also propose an instability criterion in nonlocal elasticity or dual-horizon state-based peridynamics for the fracture modeling. The criterion is formulated as the functional of nonlocal gradient in support, which minimizes the zero-energy deformation that cannot be described by the nonlocal gradient. Such an operator functional approaches zero for continuous fields but has comparable value to the strain energy density for the deformation around the crack tip. During the fracture modeling by removing particles from the neighbor list, it is safer to delete the particle with larger zero-energy deformation. The numerical examples for 2D/3D fracture modeling confirm the feasibility and robustness of this criterion. The instability criterion is possible applicable for anisotropic elastic material and hyperelastic materials.
\appendix
\section{A simple example to illustrate dual-support}\label{sec:app1}
\begin{figure}[htp]
\centering
\includegraphics[width=6cm]{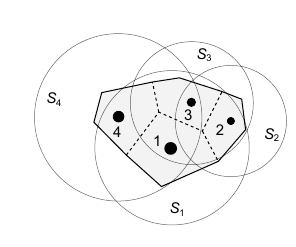}
\caption{Particles 1-4 and their supports $\cS_i, i=\{1,2,3,4\}$. }
\label{fig:Support4}
\end{figure}
In order to facilitate the comprehension of dual-support, let us consider 4 particles in Fig.\ref{fig:Support4}, each with particle volume $\Delta V_i, i=\{1,2,3,4\}$ and $\Omega=\sum_{i=1}^4 \Delta V_i$. Obviously, the support and dual-support can be listed as follows.
\begin{align}
\cS_1=\{2,3,4\}, \cS_1'=\{3,4\}\notag\\
\cS_2=\{3\}, \cS_2'=\{1,3\}\notag\\
\cS_3=\{1,2\}, \cS_3'=\{1,2,4\}\notag\\
\cS_4=\{1,3\}, \cS_4'=\{1\}\notag
\end{align}
Here we neglect whether the shape tensor is invertible or not.

The most common formula in the derivation based on NOM and variational principle is the double integrations in support and whole domain.
Consider the double integrations
\begin{align}
&\int_{\Omega} \int_{\cS_i} f_{ij} ( u_{j}- u_i) \ud V_j\ud V_i\notag\\
\approx&\sum_{i=1}^4 \Big(\sum_{j\in \cS_i} f_{ij}( u_j- u_i) \Delta V_j\Big) \Delta V_i\notag\\
=&\sum_{i=1}^4 \Big(\sum_{j\in \cS_i} f_{ij} u_j \Delta V_j\Big) \Delta V_i-\sum_{i=1}^4 \Big(\sum_{j\in \cS_i} f_{ij} u_i \Delta V_j\Big) \Delta V_i\notag
\end{align}
Expand the double summations
\begin{align}
&\sum_{i=1}^4 \Big(\sum_{j\in \cS_i} f_{ij} u_j \Delta V_j\Big) \Delta V_i\notag\\
=&\Big(f_{12} \Delta V_2 \Delta V_1 u_2+f_{13} \Delta V_3 \Delta V_1 u_3+f_{14} \Delta V_4 \Delta V_1 u_4\Big)+\Big(f_{23} \Delta V_2 \Delta V_3 u_3\Big)\notag\\
+&\Big(f_{31} \Delta V_1 \Delta V_3 u_1+f_{32} \Delta V_2 \Delta V_3 u_2\Big)
+\Big(f_{41} \Delta V_1 \Delta V_4 u_1+f_{43} \Delta V_3 \Delta V_4 u_3\Big)\notag\\
=&\Big(f_{31} \Delta V_3 +f_{41} \Delta V_4 \Big) u_1 \Delta V_1+\Big(f_{12} \Delta V_1+f_{32} \Delta V_3\Big) u_2 \Delta V_2\notag\\
+&\Big(f_{13} \Delta V_1+f_{23} \Delta V_2+f_{43} \Delta V_4\Big) u_3 \Delta V_3
+\Big(f_{14} \Delta V_1\Big) u_4 \Delta V_4\notag\\
=&\sum_{j\in S_1'} f_{j1} \Delta V_j \Delta V_1 u_1+\sum_{j\in S_2'} f_{j2} \Delta V_j \Delta V_2 u_2+\sum_{j\in S_3'} f_{j3} \Delta V_j \Delta V_3 u_3+\sum_{j\in S_4'} f_{j4} \Delta V_j \Delta V_4 u_4\notag\\
=&\sum_{i=1}^4 \Big(\sum_{j\in \cS_i'} f_{ji} \Delta V_j\Big) u_i \Delta V_i
\end{align}
Therefore
\begin{align}
&\sum_{i=1}^4 \Big(\sum_{j\in \cS_i} f_{ij} u_j \Delta V_j\Big) \Delta V_i-\sum_{i=1}^4 \Big(\sum_{j\in \cS_i} f_{ij} u_i \Delta V_j\Big) \Delta V_i\notag\\
=&\sum_{i=1}^4 \Big(\sum_{j\in \cS_i'} f_{ji} \Delta V_j\Big) u_i \Delta V_i-\sum_{i=1}^4 \Big(\sum_{j\in \cS_i} f_{ij} u_i \Delta V_j\Big) \Delta V_i\notag\\
=&\sum_{i=1}^4 \Big(\sum_{j\in \cS_i'} f_{ji} \Delta V_j-\sum_{j\in \cS_i} f_{ij} \Delta V_j\Big) u_i \Delta V_i\notag\\
\approx&\int_{\Omega} \Big(\int_{\cS_i'} f_{ji} \ud V_j-\int_{\cS_i} f_{ij} \ud V_j\Big) u_i \ud V_i
\end{align}
At last, we obtain
\begin{align}
\int_{\Omega} \int_{\cS_i} f_{ij} ( u_{j}- u_i) \ud V_j\ud V_i=\int_{\Omega} \Big(\int_{\cS_i'} f_{ji} \ud V_j-\int_{\cS_i} f_{ij} \ud V_j\Big) u_i \ud V_i
\end{align}
Above equation is widely used in the derivation of nonlocal strong form from weak form. Such expression is valid in the continuum form as well as in discrete form.
\vspace{-6pt}
\bibliographystyle{unsrt}
\bibliography{NOMfrac}
\end{document}